\documentstyle[11pt]{article}

\newcommand{\A}{ {\cal A} }
\newcommand{\B}{ {\cal B} }
\newcommand{\C}{ {\bf C} }
\newcommand{\F}{{\cal F}}

\newcommand{\R}{ {\bf R} }
\newcommand{\X}{ {\cal X} }

\newcommand{\bcop}{ \widetilde{b} }
\newcommand{\ccop}{ \widetilde{c} }
\newcommand{\xcop}{ \widetilde{x} }
\newcommand{\unucop}{ \widetilde{1} }

\newcommand{\aunun}{ a_{1}, \ldots ,a_{n} }
\newcommand{\bunun}{ b_{1}, \ldots ,b_{n} }
\newcommand{\zunun}{ z_{1}, \ldots ,z_{n} }
\newcommand{\iunuk}{ i_{1} , \ldots ,i_{k} }
\newcommand{\ziqzik}{ z_{i_{1}} \cdots z_{i_{k}} }

\newcommand{\ncps}{ ( {\cal A} , \varphi ) }
\newcommand{\noncom}{ \mbox{non-commutative probability space} }

\newcommand{\rv}{ \mbox{random variable} }
\newcommand{\rvs}{ \mbox{random variables} }
\newcommand{\ncr}{ \mbox{non-crossing} }

\newcommand{\ncpol}{ \C \langle X_{1} , \ldots ,X_{n} \rangle }
\newcommand{\coefis}{ \mbox{coef } ( \iunuk ) }
\newcommand{\seria}{ \sum_{k=1}^{\infty} \
     \sum_{i_{1}, \ldots ,i_{k} =1}^{n} \alpha_{(i_{1}, \ldots ,i_{k}) }
     z_{i_{1}}  \cdots  z_{i_{k}}   }

\newcommand{\ecpi}{ \stackrel{\pi}{\sim} }
\newcommand{\ecrho}{ \stackrel{\rho}{\sim} }
\newcommand{\eckpi}{ \stackrel{ K( \pi ) }{\sim} }
\newcommand{\egdef}{ \stackrel{def}{=} }
\newcommand{\ecdef}{ \stackrel{def}{ \Leftrightarrow } }

\newcommand{\aut}{ \stackrel{\sim}{ \longrightarrow } }
\newcommand{\lineaut}{ \stackrel{\sim}{ -- } }

\newcommand{\nn}{ \{ 1, \ldots ,n \} }
\newcommand{\kk}{ \{ 1, \ldots ,k \} }
\newcommand{\cycle}{ ( 1 \rightarrow 2 \rightarrow \cdots \rightarrow
n \rightarrow 1 ) }
\newcommand{\ovpi}{ \widetilde{\pi} }
\newcommand{\ovrho}{ \widetilde{\rho} }
\newcommand{\ovsigma}{ \widetilde{\sigma} }

\newcommand{\freestar}{ \framebox[7pt]{$\star$} }

\begin{document}

\title{\bf On the multiplication of free $n$-tuples of non-commutative
random variables}
\author{Alexandru Nica \thanks{Research done while this author was
on leave at the Fields Institute, Waterloo, and the Queen's University,
Kingston, holding a Fellowship of NSERC, Canada.}   \\
Department of Mathematics  \\
University of Michigan  \\
Ann Arbor, MI 48109-1003, USA \\
(e-mail: andu@math.lsa.umich.edu)
\and Roland Speicher \thanks{Supported by a Heisenberg Fellowship
of the DFG.} \\
Institut f\"{u}r Angewandte Mathematik \\
Universit\"{a}t Heidelberg \\
Im Neuenheimer Feld 294 \\
D-69120 Heidelberg, Germany \\
(e-mail: roland.speicher@urz.uni-heidelberg.de) }
\date{ }

\maketitle

\vspace{1in}

\begin{abstract}
Let $a_{1} , \ldots , a_{n}, b_{1} , \ldots , b_{n}$ be random variables
\setcounter{page}{0}
in some (non-commutative) probability space, such that
$\{ a_{1} , \ldots , a_{n} \}$ is free from $\{ b_{1} , \ldots , b_{n} \}$.
We show how the joint distribution of the $n$-tuple
$( a_{1} b_{1} , \ldots , a_{n} b_{n} )$ can be described in terms of the
joint distributions of $( a_{1} , \ldots , a_{n} )$ and
$( b_{1} , \ldots , b_{n} )$, by using the combinatorics of the
$n$-dimensional $R$-transform. We point out a few applications that can be
easily derived from our result, concerning the left-and-right translation
with a semicircular element (see Sections 1.6-1.10) and the compression
with a projection (see Sections 1.11-1.14) of an $n$-tuple of non-commutative
random variables. A different approach to two of these applications is
presented by Dan Voiculescu in an Appendix to the paper.
\end{abstract}

\newpage

\setlength{\baselineskip}{18pt}

{\large\bf Introduction}

$\ $

The theory of free random variables was developed in a sequence of
papers of D. Voiculescu (see \cite{VDN}, or the recent survey in \cite{V5}),
as an instrument for approaching free products of operator algebras.
Its particular aspect addressed in the present paper is the one concerning
the addition and multiplication of free random variables; as shown
by Voiculescu in \cite{V2}, \cite{V3}, a powerful method in the study of
these operations is by the use of {\em transforms} that convert them
(respectively) into addition and multiplication of complex analytic
functions - or, in an algebraic framework, of formal power series in an
indeterminate $z$. The precise definitions of these transforms (called
{\em $R$-transform} for the addition problem and {\em $S$-transform}
for the multiplication problem) will be reviewed in the Sections 1.2, 1.3
below.

In the present paper we are pursuing a combination of two ideas that have
appeared recently in the study of the $R$- and $S$-transforms.

The first idea is that the connection between the $R$- and the $S$-transform
is closer than one might suspect at first glance. A way of making this precise
was pointed out in our paper \cite{NS}, in the form of the equation
$S( \mu ) = \F ( R( \mu ))$, for $\mu$ a distribution with non-vanishing
mean, and where $\F$ is a combinatorial object with a precise significance
(``the Fourier transform for multiplicative functions on non-crossing
partitions''). A byproduct of our result in \cite{NS} is the remark that
the multiplication of free random variables can be studied
{\em directly in terms of the $R$-transform}, via an equation of the form
$R( \mu_{ab} ) = R( \mu_{a} ) \freestar R( \mu_{b} )$, where $a,b$ are
free random variables in some non-commutative probability space,
$\mu_{a} , \mu_{b} , \mu_{ab}$ stand for the distributions of $a,b$ and $ab$,
respectively, and where the operation $\freestar$ is again an object with
precise combinatorial significance, ``the convolution of multiplicative
functions on non-crossing partitions''.

The second idea, appearing in \cite{S1}, \cite{N}, is that the $R$-transform
has natural {\em multidimensional analogues}. Besides handling the addition
of free $n$-tuples ($n \geq 1$) of random variables, the multidimensional
$R$-transform has the important property that
\[
(I) \ \ \ \ \ \ \ \ \ \
[ R( \mu_{a_{1} ' , \ldots , a_{m} ' , a_{1} '' , \ldots , a_{n} '' } ) ]
( z_{1}' , \ldots , z_{m} ' , z_{1} '' , \ldots , z_{n} '' ) \ = \
\ \ \ \ \ \
\]
\[
= \ [ R( \mu_{a_{1} ' , \ldots , a_{m} ' } ) ]
( z_{1}' , \ldots , z_{m} ' )  +
[ R( \mu_{a_{1} '' , \ldots , a_{n} '' } ) ]
( z_{1} '' , \ldots , z_{n} '' ),
\]
for every family $a_{1} ' , \ldots , a_{m} ' , a_{1} '' , \ldots , a_{n} ''$
of random variables in some non-commutative probability space, such that
$\{ a_{1} ' , \ldots a_{m} ' \}$ is free from
$\{ a_{1} '' , \ldots , a_{n} '' \}$. (In Equation ($I$), by e.g.
$\mu_{ a_{1} ' , \ldots , a_{m} ' }$ we understand the joint distribution
of the variables $a_{1} ' , \ldots , a_{m} '$, and
$R( \mu_{ a_{1} ' , \ldots , a_{m} ' } )$ is a formal power series in the
$m$ non-commuting variables $z_{1} ' , \ldots , z_{m} ' .)$ This property
opens a whole new array of possibilities, in the direction of analyzing
freeness for families of random variables by the use of an
{\em $R$-transform calculus,} i.e. of a formal manipulation of power series
which lands with an equation of the type of $(I).$

$\ $

Let $a_{1} , \ldots , a_{n}, b_{1} , \ldots , b_{n}$ be random variables
in some (non-commutative) probability space, such that
$\{ a_{1} , \ldots , a_{n} \}$ is free from $\{ b_{1} , \ldots , b_{n} \}$.
In the main theorem of the present paper we show how the joint distribution
of the $n$-tuple $( a_{1} b_{1} , \ldots , a_{n} b_{n} )$ can be described
in terms of the joint distributions of $( a_{1} , \ldots , a_{n} )$ and
$( b_{1} , \ldots , b_{n} )$, {\em by using the $n$-dimensional $R$-transform.}
The formula we obtain is
\[
(II) \ \ \ \ \ \
R( \mu_{ a_{1}b_{1} , \ldots , a_{n}b_{n} } ) \ = \
R( \mu_{ \aunun } ) \freestar R( \mu_{ \bunun } ),
\]
where, of course, the crucial point is that the operation $\freestar$
can be well understood combinatorially (it extends the ``convolution of
multiplicative functions on non-crossing partitions'', mentioned above in
the 1-dimensional case).

In the support of the idea that the formula $(II)$ can be really useful
in approaching the multiplication of free $n$-tuples, we present four
applications that can be easily proved from it, on the lines of the
``$R$-transform calculus'' mentioned in the paragraph containing
Equation $(I).$

Two of these applications are concerning the left-and-right translation
of a family of random variables by a centered semicircular element which is
free from the family. The phenomenon here is that, roughly, this operation
``converts orthogonality into freeness'' (see Corollary 1.8 below). For
instance, in the important case when $e_{1}, \ldots , e_{n}$ are pairwisely
orthogonal projections in a tracial non-commutative probability space
$\ncps$, and $b \in \A$ is a centered semicircular element of radius $r$
which is free from $\{ e_{1} , \ldots , e_{n} \}$, one gets that
$\{ be_{1}b, \ldots , be_{n}b \}$ is a free family, and that $be_{k}b$ is
a Poisson element of parameters $\varphi (e_{k} )$ and $r^{2} /4$,
$1 \leq k \leq n$ (see Remark 1.9). Instead of ``left-and-right
translation with a semicircular element'' one can use in this result
``conjugation with a circular element'' or ``left-and-right translation
with a quarter-circular element'' (Remark 1.7). Moreover, we point out the
fact that if $b$ is a centered semicircular element free from the family
$\{ \aunun \}$, then, without any additional assumption on $\aunun$, one
always gets that $\{ba_{1}b, \ldots , ba_{n}b \}$ is free from
$\{ \aunun \}$ (Application 1.10).

The other two applications are related to the compression of a family of
random variables by a projection which is free from the family. The
phenomenon here is, roughly, that ``freeness is preserved by the compression''
(see Corollary 1.12 and Application 1.13 below). In the 1-dimensional case
it is interesting to note that, given a probability measure $\mu$ with
compact support on $\R$, one gets a realization of the semigroup of measures
$( \mu_{t} )_{t}$ having $R( \mu_{t} ) = tR( \mu )$ which was studied by
Bercovici and Voiculescu in \cite{BV}; namely, $\mu_{t}$ is obtained by
essentially ``compressing $\mu$ with a projection free from $\mu$ and of
trace $1/t$.'' This shows in particular that the semigroup $( \mu_{t} )_{t}$
can always be started at $t=1.$

$\ $

The paper is divided into sections as follows. In Section 1 we make a detailed
presentation of the results announced above, after reviewing the basic free
probabilistic concepts that we are using. Section 2 is devoted to some
combinatorial preliminaries about non-crossing partitions. The operation
$\freestar$ (mentioned in Eqn.$(II)$ above) is introduced and discussed in
Section 3, where the main result of the paper, Theorem 1.4, is also proved.
The four applications mentioned above (which are stated precisely in
Sections 1.6, 1.10, 1.11, 1.13) have their proofs presented in the final
Section 4.

$\ $

The work presented here was started during a workshop on operator algebra
free products and random matrices held at the Fields Institute, Waterloo,
in March 1995. We would like to thank the organizers for the very
stimulating and inspiring atmosphere that animated the workshop.

Also, we would like to acknowledge several useful discussions with
Dan Voiculescu concerning our work. It turned out that the realization
of the free Poisson variables mentioned above had been known
to him for some time (although not published); and that he could give
a different proof for the preservation of freeness among non-commutative
random variables, under compression with a free projection. His approach is
presented in an Appendix to the present paper.

$\ $

$\ $

\setcounter{section}{1}
{\large\bf 1. Presentation of the results}

$\ $

In order to make the presentation self-contained, we begin by reviewing
a few basic definitions concerning free random variables; for more details,
the reader is referred to the monograph \cite{VDN}.

$\ $

{\bf 1.1 Basic definitions} We will consider a purely algebraic framework,
where considerations on the positivity or measurability of the random
variables involved aren't necessarily required; thus by a
{\em non-commutative probability space} we will simply understand a pair
$\ncps$, where $\A$ is a complex unital algebra (``the algebra of random
variables'') and $\varphi : \A \rightarrow \C$ (``the expectation'') is
a linear functional, normalized by $\varphi (1) =1.$ The unital subalgebras
$\A_{1} , \ldots , \A_{n}$ of $\A$ are called {\em free} (with respect to
the expectation $\varphi$) if for every $k \geq 1$,
$1 \leq \iunuk \leq n$ and $a_{1} \in {\A}_{i_{1}}, \ldots ,
a_{k} \in {\A}_{i_{k}}$ such that:
\newline
(i) $i_{j-1} \neq  i_{j}$ for $1 \leq j \leq k-1$, and
\newline
(ii) $\varphi (a_{1}) = \varphi (a_{2}) = \cdots = \varphi ( a_{k} ) = 0$
\newline
it follows that $\varphi ( a_{1} a_{2} \cdots a_{k} ) = 0.$ The definition
of freeness extends to arbitrary subsets of $\A$, by defining
${\cal X}_{1} , \ldots , {\cal X}_{k} \subseteq \A$ to be free whenever
the unital algebras generated by them are so.

If $\ncps$ is a $\noncom$, and if $\aunun$ are elements (``random variables'')
in $\A$, then the {\em joint distribution} of $\aunun$ is by definition the
linear functional $\mu_{ \aunun } : \ncpol \rightarrow \C$ determined by
\begin{equation}
\left\{   \begin{array}{l}
{ \mu_{ \aunun } (1) = 1 }   \\
{ \mu_{ \aunun } ( X_{i_{1}} X_{i_{2}} \cdots X_{i_{k}} )  \ = \
\varphi ( a_{i_{1}} a_{i_{2}} \cdots a_{i_{k}} ), \mbox{ for }
k \geq 1 \mbox{ and } 1 \leq \iunuk \leq n, }
\end{array}  \right.
\end{equation}
where $\ncpol$ denotes the algebra of polynomials in $n$ non-commuting
indeterminates $X_{1}, \ldots , X_{n}$. In the case $n=1,$ the Equation (1.1)
defines the {\em distribution} of the single element $a = a_{1} \in \A$
(which is the linear functional $\mu_{a} : \C [X] \rightarrow \C$ determined
by $\mu_{a} ( X^{n} ) = \varphi ( a^{n} ),$ $n \geq 0).$

An important example of distribution is the {\em semicircle law}, which plays
in free probability a role analogous to the one of Gaussian measures in
classical probability. If $\ncps$ is a $\noncom$, an element $a \in \A$ will
be called {\em centered semicircular} of radius $r > 0$ if its distribution
is determined by
\begin{equation}
\mu_{a} ( X^{n} ) \ = \ \varphi ( a^{n} ) \ = \
\frac{2}{\pi r^2} \int_{-r}^{r} t^{n} \sqrt{r^2 - t^2} \ dt, \ \ n \geq 0.
\end{equation}

In Sections 1.7, 1.8 below we will also meet the situation of a
non-commutative probability space
$\ncps$ where $\A$ is a unital $\star$-algebra and $\varphi$ has the
property that $\varphi ( a^{*} ) = \overline{ \varphi (a) }$, for every
$a \in \A$; this is called a (non-commutative)
{\em $\star$-probability space.}
When giving the definition of a centered semicircular element $a$ in a
$\star$-probability space, one also adds to (1.2) the condition that
$a = a^{*}$.

$\ $

{\bf 1.2 The $R$-transform} ( \cite{V2}, \cite{S1}, \cite{N}) is a useful
tool for studying joint distributions of free families of random variables.
In order not to divagate too much from the main stream of our presentation,
we defer for the moment the review of the precise definition of the
$R$-transform (see Section 3.9 below), and only mention here the nature of
this object; namely, for every $n \geq 1,$ the $n$-dimensional $R$-transform
is a certain bijection $R$ from the set of linear functionals $\Sigma_n$ =
$\{ \mu : \ncpol \rightarrow \C \ | \ \mu$ linear, $\mu (1) = 1 \}$ onto
the set $\Theta_n$ of formal power series in $n$ non-commuting variables,
and without constant coefficient. (An arbitrary element of $\Theta_n$ is
thus a series of the form
\begin{equation}
f ( \zunun ) \ = \ \seria ,
\end{equation}
with $( \alpha_{(i_{1} , \ldots , i_{k})} )_{k \geq 1,
1 \leq i_{1}, \ldots , i_{k} \leq n}$ complex coefficients.)

The main property of the $R$-transform is the following: if $\ncps$ is
a non-commutative probability space,
and $a_{1} ' , \ldots a_{m} ', a_{1} '' , \ldots , a_{n} ''
\in \A$ are such that $\{ a_{1} ' , \ldots a_{m} ' \}$
\footnote[1]{ We allow the possibility that there exist $i \neq j$ such
that $a_{i} ' = a_{j} '$; $\{ a_{1} ' , \ldots a_{m} ' \}$ denotes here
the set obtained from $a_{1} ' , \ldots , a_{m} '$ after deleting the
repetitions. }
is free from
$\{ a_{1} '' , \ldots , a_{n} '' \}$, then
\begin{equation}
[ R( \mu_{a_{1} ' , \ldots , a_{m} ' , a_{1} '' , \ldots , a_{n} '' } ) ]
( z_{1}' , \ldots , z_{m} ' , z_{1} '' , \ldots , z_{n} '' ) \ = \
\end{equation}
\[
= \ [ R( \mu_{a_{1} ' , \ldots , a_{m} ' } ) ]
( z_{1}' , \ldots , z_{m} ' )  +
[ R( \mu_{a_{1} '' , \ldots , a_{n} '' } ) ]
( z_{1} '' , \ldots , z_{n} '' );
\]
and conversely, the fact that
$ R( \mu_{a_{1} ' , \ldots , a_{m} ' , a_{1} '' , \ldots , a_{n} '' } ) $
``has no cross-terms'' (i.e. it is the sum between a series in
$ z_{1}' , \ldots , z_{m} ' $  and a series in
$ z_{1} '' , \ldots , z_{n} ''$) is sufficient to ensure that
$\{ a_{1} ' , \ldots a_{m} ' \}$ is free from
$\{ a_{1} '' , \ldots , a_{n} '' \}$. Clearly, this property can make the
$R$-transform very useful for analyzing freeness, in the situations when
we have the capability of calculating it explicitly.

The $R$-transform was first remarked by Voiculescu (\cite{V1}, \cite{V2})
in the one-dimensional case, for another important property:
\begin{equation}
[ R( \mu_{a+b} ) ] (z) \ = \
[ R( \mu_{a} ) ] (z)  +  [ R( \mu_{b} ) ] (z) ,
\end{equation}
whenever $a,b$ are free in some $\noncom$ $\ncps$. This property can be
shown (\cite{S1}, \cite{N}) to hold in any dimension, i.e.:
\begin{equation}
[ R( \mu_{a_{1} + b_{1}, \ldots , a_{n} + b_{n} } ) ]
(z_{1} , \ldots , z_{n} ) \ = \
[ R( \mu_{a_{1} , \ldots , a_{n} } ) ] (z_{1} , \ldots , z_{n} ) +
[ R( \mu_{b_{1} , \ldots , b_{n} } ) ] (z_{1} , \ldots , z_{n} ) ,
\end{equation}
where $\aunun , \bunun$ are $\rvs$ in a $\noncom$ $\ncps$, and
$\{ \aunun \}$ is free from $\{ \bunun \}$. Equation (1.6) shows that, in
any case, the $R$-transform is well-suited for studying componentwise
{\em sums} of free $n$-tuples of random variables.

$\ $

{\bf 1.3 The $S$-transform} As mentioned in the introduction, our goal in the
present note is to study componentwise {\em products} of free $n$-tuples.
The case $n=1$ comes to considering the product of two free $\rvs$, and was
analyzed by Voiculescu in \cite{V3} by means of the $S$-transform. This is
a bijective map $S$ from the set of linear functionals
$\{ \mu : \C [ X ] \rightarrow \C \ | \ \mu$ linear, $\mu (1) =1,
\mu (X) \neq 0 \}$ onto the set of formal power series
$\{ g \ | \ g(z) = \sum_{n=0}^{\infty} \beta_{n} z^{n},$
$\beta_0 , \beta_1 , \beta_{2} , \ldots \in \C, \beta_0 \neq 0 \}$;
it can be described by the formula (\cite{V3}, Theorem 2.6):
\begin{equation}
[ S( \mu ) ] (z) \ = \ \frac{1+z}{z}
( \sum_{n=1}^{\infty} \mu (X^{n} ) z^n )^{ < -1 > } ,
\end{equation}
where ${ }^{ < -1 > }$ denotes the inverse under the operation of
composition of formal power series. The $S$-transform has the
``multiplicative analogue'' of the property mentioned in Eqn.(1.5), i.e.
\begin{equation}
[ S( \mu_{ab} ) ] (z) \ = \
[ S( \mu_{a} ) ] (z)  \cdot  [ S( \mu_{b} ) ] (z) ,
\end{equation}
whenever $a,b$ are free $\rvs$ in some non-commutative probability space
$\ncps$, such that $\varphi (a) \neq 0 \neq \varphi (b).$

In \cite{NS} we have found an alternative proof of (1.8), which has the
interesting feature that it goes by relating the $S$-transform to the
(1-dimensional) $R$-transform. The main steps of the argument can be
presented as follows:

(a) Consider the range-set
$\Theta_1 = \{ f \ | \ f(z) = \sum_{n=1}^{\infty} \alpha_n z^n$,
$ \alpha_1 , \alpha_2 , \ldots \in \C \}$ of the 1-dimensional
$R$-transform; we put into evidence a binary operation on $\Theta_1$,
which will be denoted here by $\freestar$, and has the property that
\begin{equation}
R( \mu_{ab} ) \ = \ R( \mu_{a} ) \freestar R( \mu_{b} )
\end{equation}
whenever $a,b$ are free $\rvs$ in a non-commutative probability space
$\ncps$ (the condition $\varphi (a) \neq 0 \neq \varphi (b)$ is not
required at this stage).

(b) We put into evidence a bijection $\F$ from
$\{ f \ | \ f(z) = \sum_{n=1}^{\infty} \alpha_n z^n$,
$ \alpha_1 , \alpha_2 , \ldots  \in \C , \alpha_1 \neq 0 \}$ onto
$\{ g \ | \ g(z) = \sum_{n=0}^{\infty} \beta_{n} z^{n},$
$\beta_0 , \beta_1 , \beta_2 , \ldots \in \C, \beta_0 \neq 0 \}$,
such that
\begin{equation}
[ \F ( f' \freestar f'' ) ] (z) \ = \
[ \F ( f') ] (z) \cdot [ \F (f'') ] (z)
\end{equation}
for every $f' , f''$ in the domain of $\F$.

(c) We show that
\begin{equation}
S( \mu_{a} ) \ = \ \F (R( \mu_{a} ))
\end{equation}
for every $\rv$ $a$ in a $\noncom$ $\ncps$, with $\varphi (a) \neq 0.$

(a),(b) and (c) above are clearly implying (1.8), since (for $a,b$ as in
(1.8)) we have
\[
S( \mu_{ab} )  \stackrel{(1.11)}{=}
\F ( R( \mu_{ab} )) \stackrel{(1.9)}{=}
\F ( R( \mu_{a} ) \freestar R( \mu_{b} )) \stackrel{(1.10)}{=}
\F ( R( \mu_{a} ))  \cdot \F ( R( \mu_{b} )) \stackrel{(1.11)}{=}
S( \mu_{a} ) \cdot S( \mu_{b} ).
\]

The framework in \cite{NS} is combinatorial, the object of study being
``the lattice of non-crossing partition of a finite ordered set'' (notion
recalled in Section 2 below). It is important to mention that the operation
$\freestar$ appearing in Eqn.(1.9) has a clear significance in this
combinatorial context, it is ``the convolution of multiplicative functions
on non-crossing partitions'' (\cite{NS}, Section 1.4).

$\ $

Now, while approaching products of free $n$-tuples of random variables via
an $n$-dimensional version of the $S$-transform seems to be difficult (or
even only partly possible - see also Remark 3.12 below), the point we would
like to make is that the operation $\freestar$ of Eqn.(1.9)
{\em does extend naturally}
to the $n$-dimensional situation, and allows us to gain information on
products of free $n$-tuples {\em via the $n$-dimensional $R$-transform.}
More precisely, we have:

$\ $

{\bf 1.4 Theorem} Let $n$ be a positive integer. There exists a binary
operation, $\freestar$, on the space of formal power series of the form
shown in (1.3), which is defined via certain ``summation formulae on
non-crossing partitions'',
\footnote[2]{ The exact definition of $\freestar$ will be given in
Section 3.2, after the necessary definitions concerning non-crossing
partitions are reviewed in Section 2.}
and has the following property: if $\ncps$ is a $\noncom$, and
$\aunun , \bunun \in \A$ are such that $\{ \aunun \}$ is free from
$\{ \bunun \}$, then
\begin{equation}
R( \mu_{ a_1 b_1 , \ldots , a_n b_n } ) \ = \
R( \mu_{ \aunun } )  \freestar  R( \mu_{ \bunun } ).
\end{equation}

$\ $

The crucial point in the above Theorem is that the operation $\freestar$
{\em has a precise combinatorial significance} (merely establishing the
existence of an operation with the property (1.12) would be trivial, but
also of no use). In the support of the idea that the result in 1.4 can be
really useful in approaching products of free $n$-tuples, we present a few
applications that can be easily derived from it.

$\ $

The first of the applications is concerning the left-and-right translation
of a family of random variables by a centered semicircular element which is
free from the family. We will use the following

$\ $

{\bf 1.5 Notation:} For $n \geq 1$ and $\mu : \ncpol \rightarrow \C$ a linear
functional normalized by $\mu (1) =1,$ we denote by $M( \mu )$ the formal
power series in $n$ non-commuting variables which has the moments of $\mu$
as coefficients, i.e.
\begin{equation}
[ M( \mu ) ] ( \zunun ) \ = \  \sum_{k=1}^{\infty} \
\sum_{i_{1} , \ldots , i_{k} =1}^{n}
\mu ( X_{i_{1}} \cdots X_{i_{k}} )
z_{i_{1}} \cdots z_{i_{k}} .
\end{equation}

$\ $

{\bf 1.6 Application} Let $\ncps$ be a $\noncom$ such that $\varphi$ is a
trace (i.e. $\varphi (xy) = \varphi (yx)$ for all $x,y \in \A$). Let
$\aunun , b \in \A$ be such that:
\newline
(i) $b$ is a centered semicircular element of radius $r$ (in the sense
reviewed in 1.1); and
\newline
(ii) $\{ \aunun \}$ is free from $b$. Then
\begin{equation}
R( \mu_{ ba_{1}b , \ldots , ba_{n}b } ) \ = \
M( \mu_{ \frac{r^2}{4} a_{1} , \ldots , \frac{r^2}{4} a_{n} } ) .
\end{equation}

$\ $

{\bf 1.7 Remark} The first step in deriving (1.14) is to note that
$\mu_{ba_{1}b, \ldots , ba_{n}b}$ =
$\mu_{ a_{1} b^{2} , \ldots ,  a_{n} b^{2} }$.
If $b$ is centered semicircular, then $b^2$ is what is called
``a Poisson element'' (see e.g. \cite{VDN}, the comment in
\footnote[3]{ A random variable in a non-commutative probability space
is called Poisson of parameters $\alpha , \beta$ if the $R$-transform of
its distribution is $\alpha  \beta z / (1 - \beta z )$;
the situation encountered here is the one
corresponding to the case $\alpha = 1, \ \beta = r^2 /4,$ with $r$ the
radius of $b$.
We mention that in the present note we have taken the liberty of multiplying
by $z$ the $R$-series as used in \cite{VDN}, because this makes the
notations easier when passing to the multidimensional case. }
(c) following to Theorem 3.7.2); hence Application 1.6 can be viewed as
concerning right (or equivalently, left) translations with this Poisson
element.

Moreover, it is known that the same Poisson element can be obtained:

- either as $q^2$, with $q$ ``quarter-circle element of radius $r$''
(i.e. having $\varphi ( q^n )$ =
\newline
$\frac{4}{\pi r^2} \int_{0}^{r} t^n \sqrt{r^2 - t^2 } \ dt$ for $n \geq 0$
- see \cite{VDN}, Definition 5.1.9);

- or as $c^{*} c$, with $c$ ``circular element of radius $r$'', i.e. of
the form $c = (x + iy)/ \sqrt{2}$ with $x,y$ free centered semicircular
elements of radius $r$ (see \cite{V4}, Definition 1.9 or \cite{VDN},
Definition 5.1.1; this makes of course sense only if
$\ncps$ is a $\star$-probability space).
\newline
Consequently, one can also replace the left-hand side of (1.14) by
$R( \mu_{qa_{1}q, \ldots ,qa_{n}q })$ or by
$R( \mu_{ca_{1}c^{*} , \ldots ,ca_{n}c^{*} })$, with $q$ and $c$ as
above.

$\ $

{\bf 1.8 Corollary} Let $\ncps$ be a non-commutative $\star$-probability
space, and let $\aunun , b$ be selfadjoint elements of $\A$ such that:

(i) $b$ is a centered semicircular element;

(ii) $a_{i} a_{j} = 0$ for $i \neq j$ (e.g., $\aunun$ can be
mutually orthogonal projections);

(iii) $\{ \aunun \}$ is free from $b$.
\newline
Then $ba_{1}b, \ldots , ba_{n}b$ form a free family in $\ncps$.

$\ $

[Alternative statements, following from Remark 1.7: one can replace
``$b$ semicircular'' by ``$q$ quarter-circular'' or by ``$c$ circular'',
where in the latter version $c$ isn't required to be selfadjoint, and
the $n$-tuple in the conclusion is $ca_{1} c^{*} , \ldots , ca_{n} c^{*}$.]

$\ $

{\bf Proof} We may assume that $\varphi$ is a trace (because  in (iii)
we are having two free Abelian subalgebras of $\A$, and by Proposition
2.5.3 of \cite{VDN}). Then
\[
[R( \mu_{ ba_{1}b , \ldots , ba_{n}b } ) ] ( \zunun ) \ = \
[M( \mu_{ \frac{r^2}{4} a_{1} , \ldots , \frac{r^2}{4} a_{n} } ) ]
( \zunun ) \ \ \mbox{ (by (1.14))}
\]
\[
= \ \sum_{m=1}^{n} ( \sum_{k=1}^{\infty}
( r^2 / 4 )^{k} \varphi (a_{m}^{k}) z_{m}^{k} )
\ \ \mbox{ (by hypothesis (ii));}
\]
hence $R( \mu_{ ba_{1}b , \ldots , ba_{n}b } )$ has no cross-terms (in
the sense reviewed in 1.2), and $ba_{1}b, \ldots , ba_{n}b$ must be free.
{\bf QED}

$\ $

{\bf 1.9 Remark} In addition to the statement of 1.8, note that the
individual distributions of the elements $ba_{1} b, \ldots , ba_{n} b$
can be effectively calculated by using the 1-dimensional case of
Eqn.(1.14); more precisely, we have
$[R( \mu_{ba_{k}b} )] (z)$ =
$[M( \mu_{a_{k}} )] ( \frac{r^2}{4} z),$ $1 \leq k \leq n,$
where $r$ is the radius of $b.$ In particular, if the moments of $a_{k}$
are given by ``a nice formula'', the same will happen with the $R$-transform
of $ba_{k} b.$ For instance, if $a_{k}$ is a projection of trace $\alpha$,
then $[M( \mu_{a_{k}} )] (z) = \alpha z / (1-z),$ hence
$[R( \mu_{ba_{k}b} )] (z) = ( \alpha \frac{r^2}{4} z)/(1- \frac{r^2}{4} z)$,
i.e. $ba_{k} b$ is a Poisson element of parameters $\alpha$ and $r^2 /4.$

$\ $

In the context of Sections 1.6-1.9, let us finally point out that we also
have:

$\ $

{\bf 1.10 Application} Let $\ncps$ be a $\noncom$ such that $\varphi$ is a
trace, and let $\aunun , b \in \A$ be such that $b$ is a centered
semicircular element, free from $\{ \aunun \}$. Then
$\{ ba_{1}b, \ldots , ba_{n}b \}$ is free from $\{ \aunun \}$.

$\ $

The other two applications we present are related to the compression by
a projection.

$\ $

{\bf 1.11 Application} Let $\ncps$ be a $\noncom$, such that $\varphi$ is a
trace. Consider an idempotent $p \in \A$, and denote
$\varphi (p) \egdef \alpha$; we assume that $\alpha \neq 0.$ If
$\aunun \in \A$ are such that $\{ \aunun \}$ is free from $p$ in $\ncps$,
then
\begin{equation}
R( \mu_{ pa_{1}p , \ldots , pa_{n}p }^{(p \A p)}  )  \ = \ \frac{1}{\alpha}
R( \mu_{ \alpha a_{1} , \ldots , \alpha a_{n} }^{( \A )}  ),
\end{equation}
where $\mu_{\ldots}^{(p \A p)}$ means that the corresponding joint
distribution is considered in the non-commutative probability space
$(p \A p , \frac{1}{\alpha} \varphi | p \A p)$
(while $\mu_{\ldots}^{( \A )}$ means that we have a joint distribution
considered in $\ncps$).

$\ $

{\bf 1.12 Corollary} Let $\ncps$ be a $\noncom$, such that $\varphi$ is
a trace, and let $p \in \A$ be an idempotent with
$\varphi (p) = \alpha \neq 0.$ If $\aunun \in \A$ is a free family of
elements of $\A$, such that $\{ \aunun \}$ is free from $p$, then
$pa_{1}p , \ldots , pa_{n}p$ is a free family in
$( p \A p, \frac{1}{\alpha} \varphi | p \A p ).$

$\ $

{\bf Proof} From (1.15) and the freeness of
$\alpha a_{1} , \ldots , \alpha a_{n}$ it is immediate that
$R( \mu_{ pa_{1}p , \ldots , pa_{n}p }^{(p \A p)}  )$ has no cross-terms;
hence $pa_{1}p , \ldots , pa_{n}p $ must be free in $p \A p.$ {\bf QED}

$\ $

The instance of the Corollary 1.12 when $\aunun$ are centered semicircular
random variables was derived in [15, Proposition 2.3] by using
arguments of random matrices; in this case, the compressions
$pa_{1} p, \ldots , pa_{n} p$ are also centered semicircular. (Note that
if an $a_{k}$, $ 1 \leq k \leq n,$ appearing in 1.12 is centered
semicircular of radius $r$ in $\ncps$, then $p a_{k} p$ is centered
semicircular of radius $r \sqrt{\alpha}$ in
$( p \A p , \frac{1}{\alpha} \varphi | p \A p );$ this is by Eqn.(1.15)
applied in the 1-dimensional case, and the known fact that the $R$-transform
of a centered semicircular element of radius $r$ is $\frac{r^2}{4} z^2$
- see \cite{VDN}, Example 3.4.4).

In \cite{V4}, Proposition 2.3 it is also remarked that if the idempotent $p$
is picked in an Abelian unital subalgebra ${\cal D} \subseteq \A$, which is
free from $\{ \aunun \}$, then $p {\cal D} p$ is free from
$\{ pa_{1} p, \ldots , pa_{n} p \}$ in $p \A p$. This latter assertion
doesn't follow from Application 1.11, but it turns out that it can
still be derived via calculations involving the operation $\freestar$ of
Theorem 1.4; moreover, the assumption that ${\cal D}$ is commutative can
be dropped - that is:

$\ $

{\bf 1.13 Application} Let $\ncps$ be a $\noncom$ such that $\varphi$ is
a trace. Let $\B \subseteq \A$ be a unital subalgebra, and let $p \in \B$
be an idempotent such that $\varphi (p) = \alpha \neq 0.$ If $\aunun \in \A$
are such that $\{ \aunun \}$ is free from $\B $ in $\ncps$, then
$\{ pa_{1}p , \ldots , pa_{n}p \}$ is free from $p \B p$ in
$(p \A p , \frac{1}{\alpha} \varphi | p \A p).$

$\ $

We conclude this section by noting that the 1-dimensional version of
Application 1.11 also has the following interesting consequence.

$\ $

{\bf 1.14 Corollary} Let $\mu$ be a compactly supported probability measure
on $\R$; then for every $t \geq 1$ there exists a unique compactly supported
probability measure $\mu_{t}$ on $\R$, such that $R( \mu_{t} ) = tR( \mu ).$
\newline
[See also \cite{BV}, Proposition 8, where the semigroup $( \mu_{t} )_{t}$
is shown to exist and exhibit strong properties, for $t$ sufficiently large.]

$\ $

{\bf Proof} It is easy to find a von Neumann algebra $\A$ with a normal
trace-state $\varphi : \A \rightarrow \C$, and $a=a^{*} \in \A$,
$(p_{\alpha} )_{0 \leq \alpha \leq 1}$ selfadjoint projections in $\A$,
such that: (i) the distribution of $a$ in $\ncps$ is $\mu$;
(ii) $\varphi (p_{\alpha} ) = \alpha$, for every $0 \leq \alpha \leq 1;$
(iii) $a$ is free from $\{ p_{\alpha} \ | \ 0 \leq \alpha \leq 1 \}$
(these elements can be realized for instance in the free product between
$L^{\infty} ( \mu )$ and the $L^{\infty}$-algebra of the Lebesgue measure
on [0,1]). For every $t \geq 1$ we consider the selfadjoint element
$a_{t} = tp_{1/t} a p_{1/t} \in p_{1/t} \A p_{1/t}$. The distribution
$\mu_{t}$ of $a_{t}$ in $p_{1/t} \A p_{1/t}$ is a compactly supported
measure on $\R$ (see \cite{VDN}, Remark 2.3.2), and has
$R( \mu_{t} ) = tR( \mu )$ by Eqn.(1.15). The uniqueness of $\mu_{t}$ is
clear, since $R( \mu_{t} ) = tR( \mu )$ determines the moments of
$\mu_{t}$. {\bf QED}

$\ $

$\ $

\setcounter{section}{2}
{\large\bf 2. Preliminaries about non-crossing partitions}

$\ $

\setcounter{equation}{0}

{\bf 2.1 Definitions} $1^{o}$ If $\pi = \{ B_{1} , \ldots , B_{k} \}$
is a partition of $\nn$ (i.e. $B_{1} , \ldots , B_{k}$ are pairwisely
disjoint, non-void sets, such that $B_{1} \cup \cdots \cup B_{k}$ =
$\nn$), then the equivalence relation on $\nn$ with equivalence classes
$B_{1} , \ldots , B_{k}$ will be denoted by $\ecpi$; the sets
$B_{1} , \ldots , B_{k}$ will be also referred to as the {\em blocks}
of $\pi$. The number of elements in the block $B_j$, $1 \leq j \leq k,$
will be denoted by $|B_{j}|$.

A partition $\pi$ of $\nn$ is called {\em non-crossing}
if for every $1 \leq i < j < k < l \leq n$ such that $i \ecpi k$
and $j \ecpi l$, it necessarily follows that
$i \ecpi j \ecpi k \ecpi l$. The set of non-crossing partitions of
$\nn$ will be denoted by $NC(n).$

For instance, all the partitions of $\nn$ with $n \leq 3$ are $\ncr$,
and the only partition of $\{ 1,2,3,4 \}$ which is not $\ncr$ is
$\{ \{ 1,3 \} , \{ 2,4 \} \}$. In general, the number of $\ncr$
partitions of $\nn$ can be shown (see e.g. \cite{K}) to be the
Catalan number $(2n)!/(n! (n+1)!)$.

$2^{0}$ For $\pi , \rho \in NC(n),$ we will write ``$\pi \leq \rho$''
if each block of $\rho$ is a union of blocks of $\pi$ (equivalently,
if we have the implication ``$i \ecpi j \Rightarrow i \ecrho j$'',
$1 \leq i,j \leq n$). Then ``$\leq$'' is a partial order relation on
$NC(n)$ (called the {\em refinement order}), and $(NC(n), \leq )$ can
in fact be shown to be a lattice.
\footnote[4]{ This means that every two elements
$\pi, \rho \in NC(n)$ have a lowest upper bound $\pi \vee \rho \in NC(n),$
and a greatest lower bound $\pi \wedge \rho \in NC(n).$ The operations
``$\vee$'' and ``$\wedge$'' will not be explicitly used in what follows;
thus, when we speak of a ``lattice isomorphism'' (or ``lattice
anti-isomorphism''), this may be taken just as an isomorphism (respectively
anti-isomorphism) for the order structure.}
We will use the notations
\begin{equation}
\left\{  \begin{array}{l}
{ 0_n \ = \ \{ \{ 1 \} , \{ 2 \} , \ldots , \{ n \} \}  }  \\
{ 1_n \ = \ \{ \ \{ 1  ,  2  , \ldots ,  n  \} \ \}  }
\end{array}   \right.
\end{equation}
for the minimal and maximal element of $NC(n)$, respectively.

The lattice of non-crossing partitions was introduced by G. Kreweras in
\cite{K}, and its combinatorics has been studied by several authors
(see e.g. \cite{SU}, and the list of references there). We will
only insist here on one basic concept (the Kreweras complementation map),
which will be extensively used in the Sections 3 and 4 of the paper.

$\ $

{\bf 2.2 The complementation map of Kreweras} is a remarkable
lattice anti-isomorphism $K:NC(n) \rightarrow NC(n)$, introduced in
$\cite{K}$, Section 3, and described as follows.

We will use the {\em circular representation} of a partition
$\pi = \{ B_{1} , \ldots , B_{k} \}$ of $\nn$, which consists in drawing
$n$ equidistant and clockwisely ordered points $P_{1} , \ldots , P_{n}$
on a circle, and in drawing, for each block $B_{j}$ ($1 \leq j \leq k$)
of $\pi$, the convex hull $H_{j}$ of the points
$\{ P_{m} \ | \ m \in B_{j} \}$. It is easily verified that $\pi$ is
$\ncr$ if and only if the $k$ convex hulls $H_{1}, \ldots , H_{k}$ are
pairwisely disjoint. Moreover, if $\pi$ is $\ncr$, then - denoting by $D$
the convex hull of the whole circle - it is also easily verified that
$D \setminus ( H_{1} \cup \cdots \cup H_{k} )$ has exactly $n+1-k$
connected components ${\widetilde{H}}_{1} , \ldots , {\widetilde{H}}_{n+1-k}$,
each of them being itself a convex set. (In order to count the connected
components of $D \setminus ( H_{1} \cup \cdots \cup H_{k} )$, one can proceed
by drawing $H_{1} , \ldots , H_{k}$ one by one. First, $D \setminus H_{1}$
has $|B_{1}|$ connected components, then the drawing of each $H_{j}$,
$2 \leq j \leq k$, increases the number of connected components of the
complement by $|B_{j}| - 1;$ so in the end the complement has
$|B_{1}| + (|B_{2}| -1) + \cdots + (|B_{k}| - 1)$ =
$(|B_{1}| + \cdots + |B_{k}|) - (k-1) =n+1-k$ components.)

Now, let $\pi = \{ B_{1} , \ldots , B_{k} \}$ be in $NC(n),$ and consider
the points $P_{1} , \ldots , P_{n}$ and the decomposition
$D = (H_1 \cup \cdots \cup H_{k}) \cup ( {\widetilde{H}}_1 \cup \cdots
\cup {\widetilde{H}}_{n+1-k} )$, as described in the preceding paragraph.
Denote the midpoints of the arcs $P_1 P_2 , \ldots , P_{n-1} P_n , \ P_n P_1$
by $Q_{1}, \ldots , Q_{n}$, respectively. Then:

$\ $

{\bf Definition} The {\em Kreweras complement} of $\pi$, denoted by
$K( \pi )$, is the partition of $\nn$ determined by
\begin{equation}
i \eckpi j \ \ecdef \ \left\{ \begin{array}{c}
{ Q_i \mbox{ and } Q_j \mbox{ belong to the same convex } }  \\
{\mbox{ set ${\widetilde{H}}_{l}$,for some $1 \leq l \leq n+1-k$ } }
\end{array}  \right\}
\end{equation}

\vspace{10pt}

$K( \pi )$ is $\ncr$ , as it is obvious by looking at its circular
representation based on the points $Q_1 , \ldots , Q_n$; note also that
$K( \pi )$ has exactly $n+1-k$ blocks (because, as it is easy to see, each
${\widetilde{H}}_{l}$, $1 \leq l \leq n+1-k,$ must contain at least one
point $Q_m$).

As a concrete example, Figure 1 illustrates that
$K( \{ \{ 1,4,8 \} , \{ 2,3 \} ,$
$\{ 5,6 \} , \{ 7 \} \} )$ =
$\{ \{ 1,3 \} , \{ 2 \} , \{ 4,6,7 \} , \{ 5 \} , \{ 8 \} \}
\in NC(8).$

\vspace{3.4in}

The following statement is a mere reformulation of the definition of the
Kreweras complement, and its proof is left as an exercise.

$\ $

{\bf Proposition} Let $\pi$ and $\rho$ be in $NC(n)$. Denote by
$\pi '$ and $\rho '$ the partitions of $\{ 2,4,6, \ldots ,$
$2n \}$ and $\{ 1,3,5, \ldots , 2n-1 \}$, respectively, which get
identified to $\pi$ and $\rho$ via the order-preserving bijections
$ \nn \rightarrow \{ 2,4,6, \ldots ,2n \}$ and
$ \nn \rightarrow \{ 1,3,5, \ldots ,2n-1 \}$. Denote by $\sigma$ the
partition of $\{ 1,2,3, \ldots , 2n \}$ formed by $\pi '$ and $\rho '$
together. Then $\sigma$ is $\ncr$ if and only if $\pi \leq K( \rho ).$

$\ $

The fact, mentioned at the beginning of this subsection, that
$K : NC(n) \rightarrow NC(n)$ is an anti-isomorphism, for every $n \geq 1,$
can be proved for instance in the following way: first, it is immediately
seen that $K^{2} ( \pi )$ is (for every $\pi \in NC(n)$) the anti-clockwise
rotation of $\pi$ with $360^{o}/n$, and this shows in particular that
$K$ is a bijection; then, it is also easy to check that
$\pi \leq \rho \Rightarrow K( \pi ) \geq K( \rho )$ - and the
converse must also hold, since $K^{2}$ is an order-preserving
isomorphism of $NC(n).$

$\ $

{\bf 2.3 A relation with the permutation group} Let ${\cal S}_n$
denote the group of all permutations of $\nn$. For
$B = \{ i_1, i_2 , \ldots , i_m \} \subseteq \nn$,
with $1 \leq i_1 < i_2 < \cdots < i_m \leq n,$ we denote by
$\tau_{B} \in {\cal S}_n$ the cycle given by
\begin{equation}
\left\{  \begin{array}{l}
{ \tau_{B} ( i_1 ) = i_2 , \ldots , \tau_{B} ( i_{m-1} ) = i_m ,
\tau_{B} ( i_{m} ) = i_1 ,}   \\
{ \tau_{B} (j) = j \mbox{ for } j \in \nn \setminus B }
\end{array}  \right.
\end{equation}
(if $|B| =1,$ then $\tau_{B}$ is of course the unit element of ${\cal S}_n$;
at the other end, $\tau_{\nn}$ is the cycle $\cycle$ ).

We will make use of a remarkable injective map
from $NC(n)$ into ${\cal S}_n$, which is denoted here by $Perm,$
and is described as follows:

$\ $

{\bf Definition} Let $\pi = \{ B_1 , \ldots , B_k \}$ be a $\ncr$ partition
of $\nn$. We denote by $Perm ( \pi ) \in {\cal S}_n$ the permutation
which has cycle decomposition:
$Perm ( \pi )  =   \tau_{B_{1}} \cdots \tau_{B_{k}}$,
with $\tau_{B_j}$, $1 \leq j \leq k,$ as defined in (2.3).

$\ $

The embedding $Perm :  NC(n) \rightarrow {\cal S}_n$ was introduced and
studied by Ph. Biane in \cite{B}, where $Perm ( \pi )$ is referred to as
``the trace of the cycle $\cycle$ on the partition $\pi$.'' The result
from $\cite{B}$ which will be needed here is that the Kreweras
complementation map has a nice interpretation in this context, namely
we have (\cite{B}, Section 1.4.2):
\begin{equation}
Perm ( \pi ) \cdot Perm ( K( \pi )) \ = \ Perm (1_n ) \  ( =  \ \cycle \  ),
\end{equation}
for every $\pi \in NC(n)$.

$\ $

{\bf 2.4 The relative Kreweras complement} Let $\rho$ be a fixed $\ncr$
partition of $\nn$. We will need a version of the Kreweras complementation
map which is considered ``relatively to $\rho$'', i.e. it is an
anti-automorphism $K_{\rho}$ of the sublattice
$\{ \pi \in NC(n) \ | \ \pi \leq \rho \}$ of $NC(n)$ (the usual Kreweras
complementation will correspond to the case $\rho = 1_n$). Roughly speaking,
the relative complement $K_{\rho} ( \pi )$ of $\pi \leq \rho$ is obtained by
looking at how $\pi$ splits the blocks of $\rho$, and then taking the usual
Kreweras complement inside each block. A more formal definition can be made
as follows.

Remark first that we have a natural lattice-isomorphism
\begin{equation}
\{ \pi \in NC(n) \ | \ \pi \leq \rho \}  \aut
NC(|B_{1}|) \times \cdots \times NC(|B_{k}|),
\end{equation}
where $B_1 , \ldots , B_k$ are the blocks of $\rho$ (and $|B_{j}|$ denotes
the number of elements of $B_{j}$). In order to describe the image of a
given $\pi \leq \rho$ under the map (2.5), one writes the blocks of $\pi$ as
$A_{1,1} , \ldots , A_{1,p_{1}} , \ldots , A_{k,1} , \ldots , A_{k,p_{k}}$,
such that $A_{j,1} \cup \cdots \cup A_{j,p_{j}} = B_{j}$ for
$1 \leq j \leq k.$ Now, if one relabels the elements
of $B_j$ (taken in increasing order) as $1,2, \ldots , |B_{j}|$, then the
partition $\{ A_{j,1} , \ldots , A_{j,p_{j}} \}$ of $B_j$ gets relabeled
into a partition $\pi_{j} \in NC( |B_{j}| ), \ 1 \leq j \leq k;$
the image of $\pi$ under
(2.5) is, by definition, the $k$-tuple $( \pi_{1} , \ldots , \pi_{k} )$.
It is immediate that (2.5) really is a bijection, and in fact a
lattice-isomorphism.

Then we can write:

$\ $

{\bf Definition} The relative Kreweras complementation map $K_{\rho}$ is
the unique map closing the square
\begin{equation}
\begin{array}{ccccc}
{ \{ \pi \in NC(n) \ | \ \pi \leq \rho \} }  &
{ -- }  &  {  \lineaut  }  &  {  - \rightarrow } &
{ NC(|B_{1}|) \times \cdots \times NC(|B_{k}|) }  \\
{ | }        &  &  &  &  { | }  \\
{ | }        &  &  &  &  { | }  \\
{ K_{\rho} } &  &  &  &  { K \times \cdots \times K }  \\
{ | }        &  &  &  &  { | }  \\
{ \downarrow }        &  &  &  &  { \downarrow }  \\
{ \{ \pi \in NC(n) \ | \ \pi \leq \rho \} }  &
{ -- }  &  {  \lineaut  }  &  { -  \rightarrow } &
{ NC(|B_{1}|) \times \cdots \times NC(|B_{k}|) }  \\
\end{array}
\end{equation}
where on the horizontals we have the isomorphism of (2.5), and on the right
vertical we have the direct product of the (usual) complementation maps on
$NC(|B_{1}|) , \ldots , NC(|B_{k}|).$

$\ $

{\bf 2.5 Relation with the permutation groups (continued)} The Equation
(2.4) mentioned in Section 2.3 is also having a relativized analogue, namely
\begin{equation}
Perm ( \pi ) \cdot Perm ( K_{\rho} ( \pi ) ) \ = \ Perm ( \rho ) ,
\end{equation}
for every $\pi , \rho \in NC(n)$ such that $\pi \leq \rho$.
The verification of (2.7) is easily done by using (2.4) and (2.6), and is
left to the reader.

The following partial converse to (2.7) (or rather the Corollary following
to it) will be used in Section 3.5.

$\ $

{\bf Proposition} Let $\ovpi , \ovsigma , \ovrho$ be non-crossing partitions
of $\nn$, such that $\ovsigma \leq \ovrho$, and
\begin{equation}
Perm ( \ovpi )  \cdot Perm ( \ovsigma ) \ = \ Perm ( \ovrho ) .
\end{equation}
Then we also have that $\ovpi \leq \ovrho$, and moreover, that
$K_{ \ovrho } ( \ovpi ) = \ovsigma$.

$\ $

{\bf Proof} The blocks of $\ovpi$ are recaptured as the orbits of $\nn$
under the action of the permutation $Perm ( \ovpi ).$ On the other hand,
the blocks of $\ovrho$ are invariant under the action of both
$Perm ( \ovrho )$ and $Perm ( \ovsigma )$ (where the latter assertion
follows from the hypothesis that $\ovsigma \leq \ovrho$); hence the
blocks of $\ovrho$ must also be invariant for
$Perm ( \ovrho ) \cdot ( Perm ( \ovsigma ))^{-1}$
$\stackrel{(2.8)}{=}  \ Perm ( \ovpi ).$ This clearly entails that every
block of $\ovrho$ is a union of blocks of $\ovpi$, i.e. that
$\ovpi \leq \ovrho$. Finally, by comparing (2.7) with (2.8) we obtain that
$Perm ( \ovsigma ) = Perm ( K_{ \ovrho } ( \ovpi )),$ hence that
$\ovsigma = K_{\ovrho} ( \ovpi ).$ {\bf QED}

$\ $

{\bf Corollary} Let $\pi , \rho$ be $\ncr$ partitions of $\nn$, such that
$\pi \leq \rho$. Then we have that $K_{\rho} ( \pi ) \leq K( \pi )$, and
moreover, that
\begin{equation}
K_{K( \pi )} ( K_{ \rho } ( \pi )) \ = \ K( \rho ) .
\end{equation}

$\ $

{\bf Proof} It suffices to verify the equality
\begin{equation}
Perm ( K_{ \rho } ( \pi )) \cdot Perm ( K( \rho )) \ = \
Perm ( K( \pi )),
\end{equation}
because we also know that $K( \rho ) \leq K( \pi )$ (as implied by the
hypothesis $\pi \leq \rho$, and the fact that $K( \cdot )$ is
order-reversing), so after that we will just have to apply the preceding
Proposition, with $\ovpi = K_{ \rho } ( \pi ),$
$\ovsigma = K( \rho ),$ $\ovrho = K( \pi ).$ But(2.10) follows immediately
if we replace $Perm (K( \rho ))$ and $Perm (K( \pi ))$ from (2.4), and
$Perm (K_{\rho} ( \pi ))$ from (2.7). {\bf QED}

$\ $

{\bf 2.6 Non-crossing pairings} A partition $\pi$ of $\{ 1,2, \ldots ,2n \}$
is called a {\em pairing} if every block of $\pi$ has exactly two elements.
We denote the set of all $\ncr$ pairings of $\{ 1,2, \ldots ,2n \}$ by
$NCP(2n).$ Non-crossing pairings have been studied for quite some time before
non-crossing partitions were introduced, see e.g. \cite{T}; it is interesting
that they are also counted by the Catalan numbers,
$|NCP(2n)| = (2n)!/(n! (n+1)!)$ for $n \geq 1.$

In particular, we have $|NC(n)| = |NCP(2n)|$ for every $n \geq 1.$ This is
not a coincidence, and the following ``canonical'' bijection between $NC(n)$
and $NCP(2n)$ will be needed in the proof of the Application 1.10.

In the next Proposition, in order to avoid any possibility of confusion,
we write $K_{(n)}$ and $K_{(2n)}$ for the Kreweras complementation maps on
$NC(n)$ and $NC(2n)$, respectively.

$\ $

{\bf Proposition} Let $n$ be a positive integer. For every $\rho \in NC(n)$,
we denote by $Twice( \rho )$ the partition of $\{ 1,2, \ldots ,2n \}$ which
is obtained by ``letting $\rho$ work on $\{ 1,3, \ldots , 2n-1 \}$ and
letting $K_{(n)} ( \rho )$ work on $\{ 2,4, \ldots , 2n \}$'', in the same
sense as in the Proposition in 2.2. Then $Twice( \rho )$ is in
$NC(2n)$, and the map
\begin{equation}
NC(n) \ni \rho \ \longrightarrow \ K_{(2n)} ( Twice ( \rho ) ) \in NC(2n)
\end{equation}
is a bijection from $NC(n)$ onto $NCP(2n).$

$\ $

{\bf Proof} $Twice ( \rho )$ is non-crossing by the Proposition in 2.2.
The map (2.11) is obviously one-to-one from $NC(n)$ into $NC(2n)$, so
we only need to verify that it takes values in $NCP(2n)$ (a one-to-one
map from $NC(n)$ into $NCP(2n)$ is necessarily onto, since
$|NC(n)| = |NCP(2n)| ).$

So, let us fix $\rho \in NC(n)$, and let us show that
$K_{(2n)} ( Twice ( \rho )) \in NCP(2n).$ The simplest way of doing this is
probably via the connection with the permutation groups discussed in 2.3,
2.5. We denote $Perm ( \rho ) \in {\cal S}_{n}$ by $\xi$, and let us also
use the notations $\gamma_{n}, \gamma_{2n}$ for the cycles
$Perm (1_{n} ) = \cycle \in {\cal S}_{n}$ and
$Perm (1_{2n} ) = (1 \rightarrow 2 \rightarrow \cdots \rightarrow 2n
\rightarrow 1) \in {\cal S}_{2n}$, respectively. Then
$Perm ( K_{(n)} ( \rho )) = \xi^{-1} \gamma_{n}$ by (2.4), and it is
immediate that $Perm (Twice ( \rho )) \in {\cal S}_{2n}$ is acting by
\begin{equation}
\left\{  \begin{array}{llll}
2k-1 & \rightarrow &  2 \xi (k) -1,     &  1 \leq k \leq n;   \\
2k   & \rightarrow &  2 \xi^{-1} (k+1), &  1 \leq k <n;    \\
2n   & \rightarrow &  2 \xi^{-1} (1) .  &
\end{array} \right.
\end{equation}

Now, $Perm(K_{(2n)} (Twice( \rho )))=Perm(Twice( \rho ))^{-1} \cdot
\gamma_{2n}$, also by (2.4). Taking (2.12) into account, we obtain that
$Perm(K_{(2n)} (Twice( \rho )))^{-1}= \gamma_{2n}^{-1} \cdot
Perm(Twice( \rho ))$ acts by:
\begin{equation}
\left\{  \begin{array}{llll}
2k-1 & \rightarrow &  2 \xi (k) -2, & \mbox{for } 1 \leq k \leq n,
\ k \neq \xi^{-1} (1);   \\
2k-1 & \rightarrow &  2n ,          & \mbox{if } k= \xi^{-1} (1); \\
2k   & \rightarrow &  2 \xi^{-1} (k+1)-1, &  1 \leq k <n;    \\
2n   & \rightarrow &  2 \xi^{-1} (1) -1.  &
\end{array} \right.
\end{equation}
As one can check without difficulty, (2.13) shows that
$Perm(K_{(2n)} (Twice( \rho )))^{-1}$ is a product of $n$ disjoint
transpositions. Hence the same holds for $Perm(K_{(2n)} (Twice( \rho )))$
itself, and this means exactly that
$K_{(2n)} ( Twice ( \rho )) \in NCP(2n).$ {\bf QED}

$\ $

$\ $

\setcounter{section}{3}
{\large\bf 3. The operation} $\freestar$ {\large\bf on formal power series}

$\ $

\setcounter{equation}{0}

{\bf 3.1 Notations} Let $n$ be a positive integer. Recall from Section 1.2
that we denote by $\Theta_{n}$ the set of formal power series without
constant coefficient in $n$ non-commuting variables $z_1 , \ldots , z_n$
(i.e., series of the form appearing in Eqn.(1.3) above). For
$f \in \Theta_{n}$ and $k \geq 1$, $1 \leq \iunuk \leq n,$ we will denote
\begin{equation}
[ \coefis ] (f) \ \egdef \  \mbox{the coefficient of $\ziqzik$ in $f$} .
\end{equation}

The following conventions for denoting coefficients will also turn out to
be handy in what follows.

(a) Let $k \geq 1$ and $1 \leq \iunuk \leq n$ be integers, and let
$B = \{ b_1 , b_2 ,  \ldots , b_l \}$ be a non-void subset of
$\{ 1, \ldots , k \}$, where $1 \leq b_1 < b_2 <  \cdots < b_l \leq k.$
Then by ``$( \iunuk ) | B$'' we will understand the $l$-tuple
$(i_{b_{1}} , i_{b_{2}} ,  \ldots , i_{b_{l}} ).$ (E.g.,
$(i_1 , i_2 , i_3 , i_4 , i_5 ) | \{ 2,3,5 \} = (i_2 , i_3 , i_5 ).)$

(b) Let $k \geq 1$ and $1 \leq \iunuk \leq n$ be integers, and let
$\pi$ be a partition of $\kk$. Then for every $f \in \Theta_{n}$ we put
\begin{equation}
[ \coefis ; \pi ] (f) \ \egdef \ \prod_{B \ block \ of \ \pi}
[ \coefis |B ] (f) .
\end{equation}
(Thus, for example, if $k=4$ and $\pi = \{ \{ 1,3 \} , \{ 2 \} , \{ 4 \} \}$,
then
\[
[ \mbox{coef } (i_1 , i_2 , i_3 , i_4 ) ; \pi ] (f) \ = \
[ \mbox{coef } (i_1 , i_3 ) ] (f) \cdot
[ \mbox{coef } (i_2 ) ] (f) \cdot
[ \mbox{coef } (i_4 ) ] (f) ,
\]
for every $f \in \Theta_{n}$ and $1 \leq i_1 , i_2 , i_3 , i_4 \leq n.)$

$\ $

Note that if $\pi$ in (3.2) is the partition (denoted in (2.1) by $1_k$)
of $\kk$ into only one block, then $[ \coefis ; \pi ] (f)$ is just the
coefficient $[ \coefis ] (f)$ appearing in (3.1). It is clear that (given
$k$ and $\iunuk$) this is the only choice of $\pi$ for which the functional
$[ \coefis ; \pi ] : \Theta_{n} \rightarrow \C$ is linear.

$\ $

{\bf 3.2 Definition} Let $n$ be a positive integer.
We denote by $\freestar$ the binary operation on the set $\Theta_{n}$ of
3.1, which is determined by the equation
\begin{equation}
[ \coefis ] (f \freestar g) \ = \  \sum_{\pi \in NC(k)}
[ \coefis ; \pi ] (f) \cdot [ \coefis ; K( \pi ) ] (g) ,
\end{equation}
holding for every $k \geq 1$ and $1 \leq \iunuk \leq n,$ and where
$K: NC(k) \rightarrow NC(k)$ is the Kreweras complementation map reviewed
in Section 2.2.

$\ $

{\bf 3.3 Remark} While (3.3) determines $f \freestar g$ completely, it is
useful to note that it can be extended to:
\begin{equation}
[ \coefis ; \rho ] (f \freestar g) \ = \
\sum_{ \begin{array}{c}
{\scriptstyle \pi  \in NC(k)} \\
{\scriptstyle \pi  \leq \rho  }
\end{array}  }
[ \coefis ; \pi  ] (f) \cdot [ \coefis ; K_{\rho } ( \pi  ) ] (g) ,
\end{equation}
holding for $k \geq 1, \ 1 \leq \iunuk \leq n,$ and some arbitrary
partition $\rho  \in NC(k).$ Eqn.(3.4) is obtained in a straightforward
way from (3.3) and the considerations in Section 2.4.

$\ $

{\bf 3.4 Remark} If in the preceding definition we take $n=1,$ then the
operation $\freestar$ (on $\Theta_{1}$) has a definite combinatorial
significance, and was analyzed in \cite{S1}, \cite{NS} in an approach where
it is called ``the convolution of multiplicative functions in the large
incidence algebra on non-crossing partitions''. It would be possible to
adapt this approach to the multidimensional situation, and place our
considerations in the framework of what is called ``the Moebius inversion
theory'' of a certain partially ordered set of ``colored $\ncr$ partitions''.
However, since the properties of $\freestar$ that are needed here can be
derived in a self-contained and elementary way, we have chosen not to enter
into any details in this direction (the reader interested in having a look
at multiplicative functions on $\ncr$ partitions is referred to \cite{S1},
Section 3, or \cite{NS}, Section 1; for the general theory of Moebius
inversion on partially ordered sets, see e.g. \cite{Sta}, Chapter 3).

$\ $

{\bf 3.5 Proposition} Let $n$ be a positive integer.
The binary operation $\freestar$ on $\Theta_{n}$ defined in 3.2 is
associative.

$\ $

{\bf Proof} Consider (and fix) the series $f,g,h \in \Theta_{n}$ and the
numbers $k \geq 1$, $1 \leq \iunuk \leq n$, about which we want to show that
\begin{equation}
[ \coefis ] ( ( f \freestar g ) \freestar h ) \ = \
[ \coefis ] ( f \freestar ( g \freestar h ) ) .
\end{equation}

We have:
\[
[ \coefis ] ( ( f \freestar g ) \freestar h )
\ \stackrel{(3.3)}{=}
\]
\[
= \  \sum_{\rho  \in NC(k)}
[ \coefis ; \rho ] (f \freestar g) \cdot [ \coefis ; K( \rho  ) ] (h)
\ \stackrel{(3.4)}{=}
\]
\begin{equation}
= \ \sum_{ \begin{array}{c}
{\scriptstyle \pi \leq \rho } \\
{\scriptstyle in \ NC(k) }
\end{array}  }
[ \coefis ; \pi ] (f) \cdot [ \coefis ; K_{\rho} ( \pi ) ] (g)
\cdot [ \coefis ; K( \rho ) ] (h).
\end{equation}

In a similar way, we see that
\[
[ \coefis ] ( f \freestar ( g \freestar h ) )
\ \stackrel{(3.3)}{=}
\]
\[
= \ \sum_{\pi \in NC(k)} [ \coefis ; \pi ] (f) \cdot
[ \coefis ; K( \pi ) ] (g \freestar h) \ =
\]
(via the substitution $K( \pi ) = \sigma$)
\[
= \ \sum_{\sigma  \in NC(k)} [ \coefis ; K^{-1} ( \sigma ) ] (f) \cdot
[ \coefis ;  \sigma  ] (g \freestar h)
\ \stackrel{(3.4)}{=}
\]
\begin{equation}
\sum_{ \begin{array}{c}
{\scriptstyle \tau  \leq \sigma } \\
{\scriptstyle in  \ NC(k)}
\end{array}  }
[ \coefis ; K^{-1} ( \sigma ) ] (f) \cdot
[ \coefis ;  \tau ] (g) \cdot [ \coefis ; K_{\sigma} ( \tau ) ] (h).
\end{equation}

We are left to establish the equality of the sums in (3.6) and (3.7).
We will do this by showing that: the map
\begin{equation}
( \pi  , \rho  ) \ \longrightarrow \ ( \tau , \sigma ) \egdef
( K_{\rho} ( \pi ) , K( \pi ) )
\end{equation}
is a bijection from $\{ ( \pi  , \rho ) \ | \ \pi , \rho  \in NC(k),
\ \pi \leq \rho \}$ onto itself, which identifies the sums (3.6) and
(3.7) term by term.

The fact that the map (3.8) really takes its domain into itself follows
from Corollary 2.5. In order to show that (3.8) is a bijection, it suffices
to check injectivity; and indeed, from
$( K_{\rho} ( \pi ) , K( \pi ) ) = ( K_{\rho '} ( \pi ' ) , K( \pi ' ) )$
we get first that $K( \pi ) = K( \pi ') \Rightarrow  \pi = \pi '$,
and then (in the notations of 2.5):
\[
K_{\rho} ( \pi ) = K_{\rho '} ( \pi ) \ \stackrel{(2.7)}{ \Rightarrow }
\ (Perm ( \pi ))^{-1} (Perm ( \rho )) = (Perm ( \pi ))^{-1} (Perm ( \rho '))
\]
\[
\Rightarrow \ Perm ( \rho ) = Perm ( \rho ') \ \Rightarrow \ \rho = \rho '.
\]
Finally, let us verify that for every $\pi \leq \rho$ in $NC(k)$, the term
in (3.6) corresponding to $( \pi , \rho )$ is equal to the term in (3.7)
corresponding to $( \tau , \sigma ) \egdef (K_{\rho} ( \pi ) , K( \pi ));$
this comes, clearly, to verifying that $\pi = K^{-1} ( \sigma )$,
$K_{\rho} ( \pi ) = \tau$, $K( \rho ) = K_{\sigma} ( \tau ).$ And indeed,
the first two of the latter equalities are obvious, while the third one
coincides with (2.9) of Corollary 2.5. {\bf QED}

$\ $

The argument proving the next proposition is very similar to the one
used to establish the Moebius inversion formula in a partially ordered
set (compare for instance to \cite{Sta}, Proposition 3.6.2); we will
skip here the details of the straightforward, but rather space-consuming
proof of the $2^{o}$``$\Leftarrow$'' part.

$\ $

{\bf 3.6 Proposition} Let $n$ be a positive integer.

$1^{o}$ The binary operation $\freestar$ on $\Theta_{n}$ has a neutral
element, which is the series
\begin{equation}
Sum ( \zunun ) \egdef z_1 + \cdots + z_n .
\end{equation}

$2^{o}$ A series $f \in \Theta_{n}$ is invertible with respect to
$\freestar$ if and only if its coefficients of degree one,
$[ \mbox{coef } (i)] (f),$ $1 \leq i \leq n,$ are all non-zero.

$\ $

{\bf Proof} $1^{o}$ is immediate.

$2^{o}$``$\Rightarrow$'': If $g$ denotes the inverse of $f$ under $\freestar$,
then for every $1 \leq i \leq n$ we have:
\[
1 = [ \mbox{coef } (i)] (Sum) \ = \
[ \mbox{coef } (i)] (f \freestar g) \ \stackrel{(3.3)}{=}
[ \mbox{coef } (i)] (f) \cdot  [ \mbox{coef } (i)] (g) ,
\]
which implies that $[ \mbox{coef } (i)] (f) \neq 0.$

$2^{o}$``$\Leftarrow$'': Say that we  want to define a family of complex
numbers, $( \beta_{( \iunuk ) } )_{k \geq 1, 1 \leq \iunuk \leq n}$,
such that the series
\[
g( \zunun ) \egdef \sum_{k=1}^{\infty} \  \sum_{\iunuk =1}^{n}
\beta_{( \iunuk )} \ziqzik
\]
has the property that $g \freestar f = Sum$.
The identification of the coefficients in the latter equality comes to an
(infinite) system of equations, and each of these equations involves
(by (3.3)) a summation over a lattice of non-crossing partitions $NC(k).$
If one separates in each such summation over $NC(k)$ the term corresponding
to the partition with only one block, $\{ \ \{ 1,2, \ldots , k \} \ \}$,
then one sees without difficulty that in fact the infinite system of
equations considered does nothing else but defining the desired coefficients
$( \beta_{( \iunuk ) } )_{k \geq 1, 1 \leq \iunuk \leq n}$,
by induction on $k$; this implies, in other words, that $f$ has a (unique)
inverse on the left under $\freestar$. The existence of an inverse on the
right is shown in a similar manner, and then, as it is well-known, the
two inverses must coincide because of the associativity of $\freestar$.
{\bf QED}

$\ $

We now turn towards giving the precise definition of the multivariable
$R$-transform, which was deferred from the Section 1.2. On the line taken
here, it is convenient to define the $R$-transform in terms of the operation
$\freestar$ of 3.2 (though, of course, this doesn't correspond to the
chronological development).

We need to introduce first the version that is appropriate, in the present
framework, for the zeta and Moebius function from the Moebius inversion
theory in partially ordered sets (compare for instance to \cite{Sta},
Section 3.7).

$\ $

{\bf 3.7 Definition} Let $n$ be a positive integer. We will call
{\em zeta power series in $n$ variables,} and denote by $Zeta$ (or
$Zeta_n$, if the precisation of $n$ is needed), the element of $\Theta_n$
given by
\begin{equation}
Zeta ( \zunun ) \ = \ \sum_{k=1}^{\infty} \ \sum_{\iunuk =1}^{n}  \ziqzik .
\end{equation}
$Zeta$ is invertible in $( \Theta_n , \freestar )$, by Proposition 3.6;
its inverse will be called the {\em Moebius series in $n$ variables,} and
will be denoted by $Moeb$ (or $Moeb_n$).

$\ $

{\bf 3.8 Remark} It is not hard to write down the Moebius series explicitly,
the formula is:
\begin{equation}
Moeb ( \zunun ) \ = \ \sum_{k=1}^{\infty} \ \sum_{\iunuk =1}^{n}
(-1)^{k+1}  \frac{ (2k-2)! }{ (k-1)! k! } \ \ziqzik
\end{equation}
(note again the occurrence of the Catalan numbers).
The shortest way for deriving (3.11) goes probably by noticing that its
verification doesn't depend on $n$, and then by invoking the literature
existent in the case $n=1,$ when $Moeb$ really is the Moebius series
associated to the lattices of non-crossing partitions (see, e.g., the
Corollary 5 in Section 3 of \cite{S1}).

$\ $

{\bf 3.9 Definition} Let $n$ be a positive integer. Let $\Theta_n$
be as above, and consider, as in Section 1.2, the set of linear
functionals $\Sigma_n = \{ \mu : \ncpol \rightarrow \C \ |$
$\mu$ linear, $\mu (1) =1 \}$. For every $\mu \in \Sigma_n$ we denote,
as in Section 1.5, by $M( \mu ) \in \Theta_n$ the formal power series
which has the moments of $\mu$ as coefficients. Then the $n$-dimensional
$R$-transform is the bijection $R: \Sigma_n \rightarrow \Theta_n$
defined by the formula
\begin{equation}
R( \mu ) \ = \ M( \mu ) \freestar Moeb, \ \ \ \mu \in \Sigma_n .
\end{equation}

$\ $

Note that an equivalent way of writing Eqn.(3.12) is
\begin{equation}
M( \mu ) \ = \ R( \mu ) \freestar Zeta, \ \ \ \mu \in \Sigma_n ;
\end{equation}
(3.13) is in some sense ``the formula for the $R^{-1}$-transform'', since
the transition from $M( \mu )$ back to $\mu$ is trivial.

Following \cite{S1}, we will also refer to the coefficients of $R( \mu )$
under the name of {\em free cumulants} of $\mu$. The above Equations
(3.13), (3.12) are in fact just a way of rewriting the version of the
Equations $( \star )$ and $( \star \star )$ in [8, Section 4]
which applies to the framework considered here.

We also mention that an alternative description of the $n$-dimensional
$R$-transform, made in terms of ``creation and annihilation operators
on the full Fock space over ${\C}^{n}$'' is presented in \cite{N};
this makes the connection between the definition given above and the
original approach of Voiculescu in \cite{V2}.

$\ $

We are only left now to present the proof of the formula (1.12) in Theorem
1.4. Before doing this, we would like to make the following important

$\ $

{\bf 3.10 Remark:} Equation (1.12) in Theorem 1.4 is equivalent to:
\begin{equation}
M( \mu_{ a_1 b_1 , \ldots , a_n b_n } ) \ = \
R( \mu_{ a_1 , \ldots ,a_n } ) \freestar M( \mu_{ b_1 , \ldots , b_n } ).
\end{equation}
Indeed, (3.14) is obtained from (1.12) by $\freestar$-operating with
$Zeta$ on the right, and the converse transition is performed by
$\freestar$-operating on the right with $Moeb$.

We take the occasion to note here that $Zeta$ and $Moeb$
{\em lie in the centre} of the semigroup $( \Theta_n , \freestar \ )$; this
is immediately verified for $Zeta$ by using the formula (3.3) and the
bijectivity of the Kreweras complementation map (and then, of course,
it also follows for $Moeb = Zeta^{-1}$). As a consequence, when
$\freestar$-operating with $Zeta$ on the right in Eqn.(1.12), we can also
associate the factor $Zeta$ (on the right-hand side) to
$R( \mu_{ a_1 , \ldots , a_n } ),$ and thus bring (1.12) to another
equivalent form,
\begin{equation}
M( \mu_{ a_1 b_1 , \ldots , a_n b_n } ) \ = \
M( \mu_{ a_1 , \ldots ,a_n } ) \freestar R( \mu_{ b_1 , \ldots , b_n } ).
\end{equation}

$\ $

{\bf 3.11 The proof of Theorem 1.4} Let $\ncps$ be a $\noncom$, and
let $\aunun , \bunun \in \A$ be such that $\{ \aunun \}$ is free from
$\{ \bunun \}$. We will prove that the equality (3.14) (equivalent, as
noted above, to (1.12) of 1.4) takes place. Thus, we fix $k \geq 1$
and $1 \leq \iunuk \leq n,$ and we will show that the coefficients of
$\ziqzik$ in $M( \mu_{ a_1 b_1 , \ldots , a_n b_n } )$ and
$R( \mu_{ \aunun } ) \freestar M( \mu_{ \bunun } )$ are equal. The line of
proof is the same as in \cite{NS}, Section 3.4 (see also \cite{S2},
Section 3.4).

We have:
\[
[ \coefis ] ( M( \mu_{ a_1 b_1 , \ldots , a_n b_n } )) \ = \
\varphi ( a_{i_1} b_{i_1} \cdots a_{i_k} b_{i_k} ) \ =
\]
\begin{equation}
= \ [ \mbox{coef } (i_1 , i_1 + n, \ldots , i_k , i_k +n) ]
( M( \mu_{ \aunun , \bunun } )),
\end{equation}
where $M( \mu_{ \aunun , \bunun } ) \in \Theta_{2n}$ will be viewed as
acting in the $2n$ variables $z_1 , \ldots , z_{2n}$, with
$z_1 , \ldots , z_n$ corresponding to the $a$'s and
$z_{n+1} , \ldots , z_{2n}$ corresponding to the $b$'s. But we know that
\[
M( \mu_{ \aunun , \bunun } )  \ \stackrel{(3.13)}{=} \
R( \mu_{ \aunun , \bunun } ) \freestar Zeta_{2n} ,
\]
so (according to Eqn.(3.3) in 3.2 and the definition of $Zeta$),
(3.16) can be continued with:
\begin{equation}
\sum_{\sigma \in NC(2k)} \
[ \mbox{coef } (i_1 , i_1 + n, \ldots , i_k , i_k +n) ; \sigma ]
\left( R( \mu_{ \aunun , \bunun } ) \right) .
\end{equation}

Now, the hypothesis that $\aunun$ is free from $\bunun$ has the consequence
that
\begin{equation}
[ R( \mu_{ \aunun , \bunun } ) ] ( z_1 , \ldots , z_n ,
z_{n+1} , \ldots , z_{2n} ) \ = \
\end{equation}
\[
[ R( \mu_{ \aunun } ) ] ( z_1 , \ldots , z_n ) +
[ R( \mu_{ \bunun } ) ] ( z_{n+1} , \ldots , z_{2n} )
\]
(see Eqn.(1.4) above); this makes clear that a partition $\sigma \in NC(2k)$
can bring a non-zero contribution to the sum (3.17) only if each block of
$\sigma$ either is contained in $\{ 1,3, \ldots , 2k-1 \}$, or is contained
in $\{ 2,4, \ldots , 2k \}$. But from the Proposition in Section 2.2 it
follows immediately that the set of partitions $\sigma \in NC(2k)$ having
the latter property is in natural bijection with the set
$\{ ( \pi , \rho ) \ | \ \pi , \rho \in NC(k), \ \pi \leq K( \rho ) \}$,
in the following way: given $\pi , \rho \in NC(k),$ with $\pi \leq K( \rho ),$
the partition $\sigma \in NC(2k)$ corresponding to $( \pi , \rho )$ is
obtained by ``letting $\pi$ work on $\{ 2, 4, \ldots , 2k \}$ and letting
$\rho$ work on $\{ 1,3, \ldots , 2k-1 \}$''. From the description of the
bijection and from (3.18) it is immediate that, whenever
$\pi , \rho , \sigma$ are as in the preceding phrase, we have
\[
[ \mbox{coef } (i_1 , i_1 + n, \ldots , i_k , i_k +n) ; \sigma ]
\left(  R( \mu_{ \aunun , \bunun } ) \right) \ =
\]
\begin{equation}
[ \coefis ; \rho ] \left( R( \mu_{ \aunun } ) \right) \cdot
[ \coefis ; \pi ] \left( R( \mu_{ \bunun } ) \right) .
\end{equation}

By putting together (3.16), (3.17), (3.19), we obtain that:
\[
[ \coefis ] ( M( \mu_{ a_1 b_1 , \ldots , a_n b_n } ) ) \ =
\]
\[
= \ \sum_{ \begin{array}{c}
{\scriptstyle \pi , \rho \ in \ NC(k) }    \\
{\scriptstyle such \ that }   \\
{\scriptstyle \pi \leq K( \rho ) }
\end{array}  }
[ \coefis ; \rho ] ( R( \mu_{ \aunun } )) \cdot
[ \coefis ; \pi ] ( R( \mu_{ \bunun } ))
\]
\begin{equation}
= \ \sum_{ \rho \in NC(k) } [ \coefis ; \rho ] ( R( \mu_{ \aunun } ))
\cdot (  \sum_{  \begin{array}{c}
{\scriptstyle \pi \leq K( \rho )}  \\
{\scriptstyle in \ NC(k) }
\end{array}  }
[ \coefis ; \pi ] ( R( \mu_{ \bunun } )) \ ).
\end{equation}
Finally, by using Eqn.(3.4) we infer that the second sum in (3.20) is equal to
\newline
$[ \coefis ; K( \rho ) ] ( R( \mu_{ \bunun } ) \freestar Zeta )$; and
since we know that $R( \mu_{ \bunun } ) \freestar Zeta$ is the same thing
as $M( \mu_{ \bunun } ),$ we can conclude that the quantity in (3.20) equals
\[
\sum_{\rho \in NC(k)} [ \coefis ; \rho ] ( R( \mu_{ \aunun } )) \cdot
[ \coefis ; K( \rho ) ] ( M( \mu_{ \bunun } ))
\]
\[
\stackrel{(3.3)}{=} \ [ \coefis ] \left( R( \mu_{ \aunun } ) \freestar
M( \mu_{ \bunun } ) \right) ,
\]
as desired. {\bf QED}

$\ $

{\bf 3.12 Remark} As we recalled in Section 1.3, the 1-dimensional instance
of Theorem 1.4 is related to the $S$-transform of Voiculescu, and it is
natural to ask whether (and to what extent) could the $S$-transform itself
be adapted to work in the multidimensional case. We are not able, at present,
to settle this question in a satisfactory way. Let us note that the existence
of an $n$-dimensional $S$-transform would be equivalent to the existence
of an $n$-dimensional ``combinatorial Fourier transform on $\ncr$ partitions'',
the analogue of the map $\F$ appearing in (b) and (c) of Section 1.3. The
$n$-dimensional version of $\F$ should be an isomorphism from, say,
$\{ f \in \Theta_n \ | \ [ \mbox{coef } (i) ] (f) = 1, \ 1 \leq i \leq n \}$,
endowed with $\freestar$, onto some semigroup ${\cal M}$, where the operation
on ${\cal M}$ should be in some way related to the multiplication of formal
power series. The various candidates we have tested  for $\F$ and
${\cal M}$ have all failed to be homomorphic. While, of course, it is not
impossible that a more fortunate choice can be found, we would like to remark
that some substantial difference between the 1- and multi-dimensional cases
is to be expected anyway, in view of the fact that $( \Theta_n , \freestar \ )$
is {\em commutative} if and only if $n=1.$ This distinction is particularly
significant in the light of the original approach of Voiculescu in \cite{V3},
where the $S$-transform is found as the exponential of a certain commutative
Lie group, the group-operation on which is related to $\freestar$ on
$\Theta_1$. (The analogue of this Lie group in the multidimensional case is
non-commutative, and it seems unlikely that a homomorphic $S$-transform could
be found by the same method.)

Let us also mention that in a recent work (\cite{H}, Part B), U. Haagerup
has found another proof of the multiplicativity of the $S$-transform, based
on an elegant adaptation of the ``Fock space model'' for non-commutative
random variables. It is possible that the right concept of multidimensional
$S$-transform might be found via a better understanding of this adapted model.

$\ $

$\ $

$\ $

\setcounter{section}{4}
{\large\bf 4. Proof of the Applications 1.6, 1.10, 1.11, 1.13}

$\ $

\setcounter{equation}{0}

{\bf 4.1 Notation} Let $n$ be a positive integer. We consider the set
$\Theta_n$ of formal power series without constant coefficient in $n$
non-commuting variables $z_1 , \ldots , z_n$, and use the notations for
coefficients introduced in Section 3.1.

For $f \in \Theta_n$ and a number $r \neq 0$ we will denote by
$f \circ D_r$, and call {\em the dilation of f by r,} the series in
$\Theta_n$ defined by
\[
(f \circ D_r ) ( z_1 , \ldots , z_n ) \ = \
f( r z_1 , \ldots , r z_n ),
\]
or equivalently, by
\begin{equation}
[ \coefis ] ( f \circ D_r ) \ = \ r^k [ \coefis ] (f),
\end{equation}
for every $k \geq 1$ and $1 \leq \iunuk \leq n.$

It is immediate that (4.1) can be extended to
\begin{equation}
[ \coefis ; \pi ] ( f \circ D_r ) \ = \ r^k [ \coefis ; \pi ] (f),
\end{equation}
for every $k \geq 1, 1 \leq \iunuk \leq n$ and $\pi \in NC(k).$ Moreover,
from (4.2) and (3.3) it clearly follows that the operation $\freestar$
discussed in the preceding section behaves nicely under dilations, i.e.
\begin{equation}
( f \circ D_r ) \freestar g \ = \ f \freestar (g \circ D_r ) \ = \
( f \freestar g ) \circ D_r ,
\end{equation}
for every $f,g \in \Theta_n$ and $r \neq 0.$

$\ $

The next lemma will be used in the proofs of 1.6 and 1.10.

$\ $

{\bf 4.2 Lemma} Given a series $f(z) = \sum_{k=1}^{\infty} \alpha_k z^k$
in $\Theta_1$ and a positive integer $n,$ let us denote by
$f(z_1 + \cdots + z_n )$ the series in $\Theta_n$ determined by
\[
[ \coefis ] (f( z_1 + \cdots + z_n )) \ = \
\mbox{the coefficient of $z^k$ in $f$,}
\]
for every $k \geq 1$ and $1 \leq \iunuk \leq n.$ Then we have that
\begin{equation}
[ R( \mu_{ \underbrace{a, \ldots ,a}_{n} } ) ] ( \zunun ) \ = \
[ R( \mu_{a} ) ] ( z_1 + \cdots + z_n ),
\end{equation}
for every $\rv$ $a$ in a $\noncom$ $\ncps$, and for every $n \geq 1.$

Particular cases: if $b$ is a centered semicircular element of radius $r$
in the noncommutative probability space $\ncps$, then:
\begin{equation}
[ R( \mu_{ \underbrace{b, \ldots ,b}_{n} } ) ] ( \zunun ) \ = \
\frac{r^{2}}{4} \sum_{i,j =1}^{n} z_{i} z_{j} \in \Theta_n ;
\end{equation}
\begin{equation}
R( \mu_{ \underbrace{b^2 , \ldots , b^2}_{n} } ) \ = \
Zeta_n \circ D_{r^{2} /4} \in \Theta_n .
\end{equation}

$\ $

{\bf Proof} The analogue of Eqn.(4.4) with the $R$-series replaced by
the $M$-series is immediate. It is also immediately checked that
\begin{equation}
(f( z_1 + \cdots + z_n )) \freestar (g( z_1 + \cdots + z_n )) \ = \
(f \freestar g) ( z_1 + \cdots + z_n )
\end{equation}
for every $f,g \in \Theta_1$, where the $\freestar$-operations in the
left-hand and right-hand side of (4.7) are considered in $\Theta_n$ and
$\Theta_1$, respectively. Since from Eqn.(3.11) in 3.8 it is clear that
$Moeb_n ( \zunun ) = Moeb_1 ( z_1 + \cdots + z_n ),$ we have:
\[
[ R( \mu_{ a, \ldots ,a } ) ] ( \zunun ) \ = \
( M( \mu_{ a, \ldots ,a } ) \freestar Moeb_n  ) ( \zunun ) \ = \
\]
\[
= \ [ M( \mu_{a} ) ] ( z_1 + \cdots + z_n ) \freestar
Moeb_1  ( z_1 + \cdots + z_n )
\]
\[
\stackrel{(4.7)}{=}
( M( \mu_{a} ) \freestar Moeb_1 ) ( z_1 + \cdots + z_n ) \ = \
[ R( \mu_{a} ) ] ( z_1 + \cdots + z_n ).
\]

Now, if $b$ is centered semicircular of radius $r$ in $\ncps$, then
$[ R( \mu_{b} ) ] (z) = \frac{r^2}{4} z^2$ (see \cite{VDN}, Example 3.4.4);
the formula (4.5) follows from this and (4.4). Similarly, in order to prove
(4.6) we only need to show that
$[ R( \mu_{b^2} ) ] (z) = \sum_{k=1}^{\infty} ( \frac{r^2}{4} z)^{k} .$
This is equivalent to
\begin{equation}
R( \mu_{b^2} ) \freestar Zeta_{1} \ = \
( \ \sum_{k=1}^{\infty} ( \frac{r^2}{4} z)^{k} \ ) \freestar Zeta_{1} .
\end{equation}
By (3.13), the left-hand side of (4.8) is $M( \mu_{b^2} ),$ i.e. the series
\[
\sum_{k=1}^{\infty} \varphi (b^{2k} ) z^k \ = \
\sum_{k=1}^{\infty} ( \frac{2}{\pi r^2} \int_{-r}^{r} t^{2k}
\sqrt{r^2 -t^2} \ dt ) z^k .
\]
A direct calculation, which takes into account that
$|NC(k)| = (2k)! / ( k! (k+1)!),$ finds the right-hand side of (4.8) equal
to $\sum_{k=1}^{\infty} (r^2 /4)^k \cdot (2k)!/ (k! (k+1)!) \cdot z^k ;$
so (4.8) comes to proving that
\[
\frac{2}{\pi r^2} \int_{-r}^{r} t^{2k} \sqrt{r^2 -t^2} \ dt  \ = \
{ \left( \frac{r^2}{4} \right)}^k \cdot \frac{(2k)!}{k!(k+1)!} , \ \ \
k \geq 1,
\]
which is easily done via integration by parts and induction on $k.$
{\bf QED}

$\ $

{\bf 4.3 Proof of Application 1.6} Let $\ncps$ and $\aunun , b \in \A$ be as
in the statement of 1.6. By using the trace-property of $\varphi$, it is
immediately seen that
$\mu_{ ba_{1}b, \ldots , ba_{n}b } = \mu_{ a_{1}b^{2} , \ldots , a_{n}b^{2} }$.
Thus:
\[
R( \mu_{ ba_{1}b, \ldots , ba_{n}b } ) \ = \
R( \mu_{ a_{1}b^{2} , \ldots , a_{n}b^{2} }  )
\ \stackrel{(1.12)}{=} \
R( \mu_{ \aunun } ) \freestar R( \mu_{ b^2 , \ldots , b^2 }  )
\]
\[
\stackrel{(4.6)}{=} \
R( \mu_{ \aunun } ) \freestar ( Zeta \circ D_{r^{2} /4} )
\ \stackrel{(4.3)}{=} \
( R( \mu_{ \aunun } ) \freestar Zeta )  \circ D_{r^{2} /4}
\ \stackrel{(3.13)}{=} \
M( \mu_{ \aunun } ) \circ D_{r^{2} /4}  .
\]
The fact that the latter series can be also written as
$M( \mu_{ \frac{r^2}{4} a_{1} , \ldots , \frac{r^2}{4} a_{n} } )$
is immediate. {\bf QED}

$\ $

We now head towards the proof of 1.11.
We will also need to use the behavior of the operation $\freestar$ under
multiplication by a scalar, as described in the following

$\ $

{\bf 4.4 Lemma:} Let $n$ be a positive integer, and let $f,g$ be two series
in $\Theta_n$. For every $r \neq 0$ we have
\begin{equation}
(rf) \freestar (rg) \ = \ ( r (f \freestar g) ) \circ D_r .
\end{equation}

$\ $

{\bf Proof} For every $k \geq 1$ and $1 \leq \iunuk \leq n$ we have:
\[
[ \coefis ] ( (rf) \freestar (rg) )
\]
\[
\stackrel{(3.3)}{=} \
\sum_{ \pi \in NC(k) } [ \coefis ; \pi ] (rf) \cdot
[ \coefis ; K( \pi ) ] (rg)
\]
\[
\stackrel{(3.2)}{=} \
\sum_{ \pi \in NC(k) } r^{| \pi |} [ \coefis ; \pi ] (f) \cdot
r^{|K( \pi )|} [ \coefis ; K( \pi ) ] (g)
\]
\[
= \ r^{k+1} \sum_{ \pi \in NC(k) } [ \coefis ; \pi ] (f) \cdot
[ \coefis ; K( \pi ) ] (g)
\]
(because $| \pi | + | K( \pi ) | = k+1$ for every $\pi \in NC(k)$)
\begin{equation}
= \ r^{k+1} [ \coefis ] ( f \freestar g ).
\end{equation}
But it is clear that the quantity in (4.10) is exactly the coefficient
of $\ziqzik$ in $(r(f \freestar g)) \circ D_r$. {\bf QED}

$\ $

{\bf 4.5 Remark} It is convenient to use Eqn.(4.9) in the slightly
modified form:
\begin{equation}
\frac{1}{r} ( f \freestar (rg)) \ = \
(( \frac{1}{r} f) \circ D_{r} ) \freestar g, \ \ f,g \in \Theta_n, \ r \neq 0;
\end{equation}
(4.11) reduces to (4.9) via the substitution $f = r \widetilde{f}$,
$g = \widetilde{g}$.

$\ $

{\bf 4.6 Proof of Application 1.11} Let $\ncps$ and $\aunun , p \in \A$
be as in the statement of 1.11. Recall that $\alpha \neq 0$ denotes the
trace of the idempotent $p$.

It is clear that
\[
M( \mu_{ pa_{1}p, \ldots , pa_{n}p }^{(p \A p)} )  \ = \
\frac{1}{\alpha}
M( \mu_{ pa_{1}p, \ldots , pa_{n}p }^{( \A )} )
\]
(since the trace on $p \A p$ is just the restriction of the one on $\A$,
normalized by a factor of $1/ \alpha$). Also, from the trace property of
$\varphi$ it follows immediately that
$\mu_{ pa_{1}p, \ldots , pa_{n}p }^{( \A )}$ =
$\mu_{ a_{1}p, \ldots , a_{n}p }^{( \A )}$. Hence, we can write:
\[
M( \mu_{ pa_{1}p, \ldots , pa_{n}p }^{(p \A p)} ) \ = \
\frac{1}{\alpha} M( \mu_{ a_{1}p, \ldots , a_{n}p }^{( \A )} )
\ \stackrel{(3.14)}{=} \ \frac{1}{\alpha}  \left(
R( \mu_{ \aunun }^{( \A )} )  \freestar
M( \mu_{ p, \ldots ,p }^{( \A )} ) \right)
\]
\[
= \ \frac{1}{\alpha} \left(  R( \mu_{ \aunun } ) \freestar
( \alpha  Zeta )  \right)
\ \stackrel{(4.11)}{=} \ (( \frac{1}{\alpha} R( \mu_{ \aunun }^{( \A )} ))
\circ D_{\alpha} ) \freestar Zeta.
\]
But from $M( \mu_{ pa_{1}p, \ldots , pa_{n}p }^{(p \A p)} )$ =
$(( \frac{1}{\alpha} R( \mu_{ \aunun }^{( \A )} )) \circ D_{\alpha} )
\freestar Zeta$ we obtain, by $\freestar$-operating with $Moeb$ on the right:
\begin{equation}
R( \mu_{ pa_{1}p, \ldots , pa_{n}p }^{(p \A p)} ) \ = \
M( \mu_{ pa_{1}p, \ldots , pa_{n}p }^{(p \A p)} ) \freestar Moeb \ = \
( \frac{1}{\alpha} R( \mu_{ \aunun }^{( \A )} )) \circ D_{\alpha} .
\end{equation}
The fact that the rightmost expression in (4.12) can be also written as
$\frac{1}{\alpha} R( \mu_{ \alpha a_{1} , \ldots , \alpha a_{n} }^{( \A )} )$
is an easy exercise, left to the reader. {\bf QED}

$\ $

For the proof of the remaining Applications 1.10 and 1.13, we will use
the following freeness criterion.

$\ $

{\bf 4.7 Proposition} Let $\ncps$ be a $\noncom$ such that $\varphi$ is
a trace, let $\B \subseteq \A$ be a unital subalgebra, and let
$\X \subseteq \A$ be a non-void subset. Assume that for every $m \geq 1$
and $b_1 , \ldots , b_m \in \B$, $x_1 , \ldots , x_m \in \X$, it is true
that
\begin{equation}
\varphi ( b_1 x_1 b_2 x_2 \cdots b_m x_m ) \ =
\ [ \mbox{coef } (1,2, \ldots , m) ]
(R( \mu_{ b_1 , \ldots , b_m } ) \freestar M( \mu_{ x_1 , \ldots , x_m }))
\end{equation}
(where, according to the notations set in 3.1, the right-hand side of (4.13)
denotes the coefficient of $z_1 z_2 \cdots z_m$ in the formal power series
$R( \mu_{ b_1 , \ldots , b_m } ) \freestar M( \mu_{ x_1 , \ldots , x_m })
\in \Theta_m ).$ Then $\X$ is free from $\B$ in $\ncps$.

$\ $

Note that the condition expressed in (4.13) is also necessary for the freeness
of $\X$ from $\B$, as it clearly follows from Theorem 1.4 and Remark 3.10
(and the obvious fact that $ \varphi ( b_1 x_1 b_2 x_2 \cdots b_m x_m )$ =
$[ \mbox{coef } (1,2, \ldots , m)](M( \mu_{b_1 x_1 , \ldots , b_m x_m} )) \ ).$

$\ $

{\bf Proof} Let ${\cal C}$ be the unital subalgebra of $\A$ generated by
$\X$. The fact that $\X$ is free from $\B$ means, by definition, that
the subalgebras $\B$ and ${\cal C}$ are free in $\ncps$.
We consider the free product of unital algebras $\B \star {\cal C}$
(its construction is independent of the fact that $\B$ and ${\cal C}$ lie in
the same algebra $\A$); we will denote the canonical embeddings of $\B$ and
${\cal C}$ in $\B \star {\cal C}$ by $b \rightarrow \bcop$ and
$c \rightarrow \ccop$, respectively (that is, $\bcop$ ``is the name'' of the
element $b \in \B$, when viewed in $\B \star {\cal C}$, and similarly for
$c \in {\cal C}$ and $\ccop \in \B \star {\cal C} ).$ By the
universality of $\B \star {\cal C}$, there exists a unique homomorphism
of unital algebras $\Phi : \B \star {\cal C} \rightarrow \A$, such that
$\Phi ( \bcop ) =b$ for every $b \in \B$ and $\Phi ( \ccop ) =c$
for every $c \in {\cal C}$ (see e.g. \cite{VDN}, Section 1.2).

Now, let us consider on $\B \star {\cal C}$ the free product functional
$\psi = ( \varphi | \B ) \star ( \varphi | {\cal C} );$ this is a trace,
because $\varphi | \B$ and $\varphi | {\cal C}$ are so and by Proposition
2.5.3 of \cite{VDN}. We remark that:

(a) We have
\begin{equation}
\varphi ( b_1 x_1 \cdots b_m x_m ) \ = \
\psi ( \bcop_1 \xcop_1 \cdots \bcop_m \xcop_m ),
\end{equation}
for every $m \geq 1, \ b_1 , \ldots b_m \in \B , \ x_1 , \ldots ,x_m \in \X .$
(4.14) holds because both its sides are equal to:
\begin{equation}
[ \mbox{coef } (1,2, \ldots ,m) ] ( R( \mu_{ b_1 , \ldots , b_m } )
\freestar M( \mu_{ x_1 , \ldots , x_m } )) .
\end{equation}
The equality  between (4.15) and the left-hand side of (4.14) is ensured by
the hypothesis, while the equality between (4.15) and the right-hand side of
(4.14) comes out from Theorem 1.4 in the form presented in Remark 3.10 (by
also taking into account that $\{ \bcop \ | \ b \in \B \}$ and
$\{ \ccop \ | \ c \in {\cal C} \}$ are free in
$( \B \star {\cal C} , \psi )$, and the obvious fact that
$\mu_{ \bcop_1 , \ldots , \bcop_m } = \mu_{ b_1 , \ldots , b_m }$,
$\mu_{ \xcop_1 , \ldots , \xcop_m } = \mu_{ x_1 , \ldots , x_m }$).

(b) (4.14) implies that $\varphi \circ \Phi = \psi$. Indeed, since
$\Phi ( \bcop_1 \xcop_1 \cdots \bcop_m \xcop_m ) = b_1 x_1 \cdots b_m x_m$,
(4.14) is verifying the coincidence of the functionals $\varphi \circ \Phi$
and $\psi$ on elements of the form $\bcop_1 \xcop_1 \cdots \bcop_m \xcop_m$.
Since $\varphi \circ \Phi$ and $\psi$ are traces, their coincidence also
follows on elements of the form
$\bcop_1 \xcop_1 \cdots \bcop_m \xcop_m \bcop_{m+1}$, with $m \geq 0$
(if $m=0$ this is clear, if $m \geq 1$ we send
$\bcop_{m+1}$ to the front and multiply it with $\bcop_{1}$).
But the elements of these two forms are spanning linearly together all of
$\B \star {\cal C}$ (first, it is clear that $\B \star {\cal C}$ is
spanned linearly by products of $\bcop$'s and $\xcop$'s; then any such product
can be turned into an alternating one, by multiplying together the consecutive
$b$'s and by inserting an $\unucop$, $1 \in \B$, between consecutive
$\xcop$'s; also, by adding if necessary an $\unucop$ on the left, it may be
always assumed that the monomial starts with a $\bcop$).

(c) The freeness of $\B$ and ${\cal C}$ in $\ncps$ is an immediate
consequence of the fact that $\varphi \circ \Phi = \psi .$ Indeed, let us
consider an alternating sequence of elements from $\B$ and ${\cal C}$, e.g.
$c_1 , b_1, \ldots , c_n , b_n , c_{n+1}$, such that
$\varphi (c_1 ) = \varphi (b_1 ) = \cdots = \varphi (c_n ) =
\varphi (b_n ) = \varphi (c_{n+1} ) =0.$ Then, clearly, we also have
$\psi ( \ccop_1 ) = \varphi ( c_1 ) =0,$
$\psi ( \bcop_1 ) = \varphi ( b_1 ) =0,$ $\ldots$,
$\psi ( \ccop_{n+1} ) = \varphi ( c_{n+1} ) =0.$ From the freeness of
$\{ \bcop \ | \ b \in \B \}$ and $\{ \ccop \ | \ c \in {\cal C} \}$
in $( \B \star {\cal C} , \psi )$ it follows that
$\psi ( \ccop_1 \bcop_1 \cdots \ccop_n \bcop_n \ccop_{n+1} ) =0;$ but
``$\varphi \circ \Phi = \psi$'' implies that
$\psi ( \ccop_1 \bcop_1 \cdots \ccop_n \bcop_n \ccop_{n+1} )$ =
$\varphi ( c_1 b_1 \cdots c_n b_n c_{n+1} ),$ so the latter quantity is
also vanishing. {\bf QED}

$\ $

{\bf 4.8 Remark} The preceding proposition could be equivalently stated
by replacing (4.13) with
\begin{equation}
\varphi ( b_1 x_1 b_2 x_2 \cdots b_m x_m ) \ =
\ [ \mbox{coef } (1,2, \ldots , m) ]
(M( \mu_{ b_1 , \ldots , b_m } ) \freestar R( \mu_{ x_1 , \ldots , x_m })).
\end{equation}
The proof would be exactly the same, with the only detail that at the point
where Remark 3.10 is invoked (in part (a) of the proof), we would now refer
to Eqn.(3.15) instead of (3.14).

$\ $

{\bf 4.9 The proof of Application 1.13} Let $\ncps$, $\aunun \in \A$,
$p \in \B \subseteq \A$ be as in the statement of 1.13.
We will use the criterion established in 4.7 for proving the freeness of
$\X = \{ pa_{1}p, \ldots , pa_{n}p \}$ from the subalgebra $p \B p$ in
the $\noncom$ $( p \A p , \frac{1}{\alpha} \varphi | p \A p )$.
We fix $m \geq 1$ and $pb_{1}p , \ldots , pb_{m}p \in p \B p$ (with
$b_{1} , \ldots , b_{m} \in \B$), $x_1 , \ldots , x_m \in \X$, about
which we want to verify the equality (4.13) of 4.7. Writing explicitly
$x_1 = pa_{i_{1}}p, \ldots , x_m = pa_{i_{m}}p$ (for some
$ 1 \leq i_1 , \ldots i_m \leq n$), the desired equality comes to:
\begin{equation}
\frac{1}{\alpha} \varphi ( (pb_1 p) (pa_{i_{1}}p) \cdots
(pb_{m}p) (pa_{i_{m}}p)) \ = \
\end{equation}
\[
[ \mbox{coef } (1,2, \ldots ,m)]
\left( R( \mu_{ pb_{1}p, \ldots , pb_{m}p }^{(p \A p)} ) \freestar
M( \mu_{ pa_{i_{1}}p, \ldots , pa_{i_{m}}p }^{(p \A p)} ) \right) .
\]

The two power series involved in the right-hand side of (4.17) can be
re-written as follows:
\begin{equation}
R( \mu_{ pb_{1}p , \ldots , pb_{m}p }^{(p \A p)} ) \ = \
M( \mu_{ pb_{1}p , \ldots , pb_{m}p }^{(p \A p)} ) \freestar Moeb \ =
\ \left(  \frac{1}{\alpha}
M( \mu_{ pb_{1}p , \ldots , pb_{m}p }^{( \A )} ) \right)  \freestar Moeb,
\end{equation}
and
\[
M( \mu_{ pa_{i_{1}}p, \ldots , pa_{i_{m}}p }^{(p \A p)} ) \ = \
\frac{1}{\alpha}
M( \mu_{ pa_{i_{1}}p, \ldots , pa_{i_{m}}p }^{( \A )} ) \ = \
\frac{1}{\alpha} M( \mu_{ a_{i_{1}}p, \ldots , a_{i_{m}}p }^{( \A )} )
\]
\[
= \ \frac{1}{\alpha} \left(
R( \mu_{ a_{i_{1}} , \ldots , a_{i_{m}} }^{( \A )} )
\freestar  M( \mu_{ p, \ldots ,p }^{( \A )} ) \right)  \ \
\mbox{  (by Theorem 1.4 and Remark 3.10)}
\]
\begin{equation}
= \ \frac{1}{\alpha}  \left(
R( \mu_{ a_{i_{1}} , \ldots , a_{i_{m}} }^{( \A )} )
\freestar  ( \alpha Zeta )  \right) \ \stackrel{(4.11)}{=} \
\left( ( \frac{1}{\alpha}
R( \mu_{ a_{i_{1}} , \ldots , a_{i_{m}} }^{( \A )} ))
\circ D_{\alpha} \right)  \freestar Zeta.
\end{equation}
When $\freestar$-multiplying together the rightmost expressions in
(4.18) and (4.19), the factors $Moeb$ and $Zeta$ cancel (because $Moeb$
and $Zeta$ are central, and inverse to each other), and we obtain:
\[
\left(  \frac{1}{\alpha}
M( \mu_{ pb_{1}p , \ldots , pb_{m}p }^{( \A )} ) \right)  \freestar
\left(  ( \frac{1}{\alpha}
R( \mu_{ a_{i_{1}}  , \ldots , a_{i_{m}} }^{( \A )} )) \circ D_{\alpha}
\right)
\]
\[
\stackrel{(4.3)}{=} \  \left(  ( \frac{1}{\alpha}
M( \mu_{ pb_{1}p , \ldots , pb_{m}p }^{( \A )} ) ) \freestar
( \frac{1}{\alpha}
R( \mu_{ a_{i_{1}} , \ldots , a_{i_{m}} }^{( \A )} ))  \right)
\circ D_{\alpha}
\]
\[
\stackrel{(4.9)}{=} \  \left(  \frac{1}{\alpha}
( \ M( \mu_{ pb_{1}p , \ldots , pb_{m}p }^{( \A )} ) \freestar
R( \mu_{ a_{i_{1}} , \ldots , a_{i_{m}} }^{( \A )} ) \ )
\circ D_{1/ \alpha} \right)  \circ D_{\alpha}
\]
\[
= \ \frac{1}{\alpha} \left(
M( \mu_{ pb_{1}p , \ldots , pb_{m}p }^{( \A )} ) \freestar
R( \mu_{ a_{i_{1}} , \ldots , a_{i_{m}} }^{( \A )} )  \right)
\]
\[
= \ \frac{1}{\alpha} M( \mu_{ (pb_{1}p)a_{i_{1}} , \ldots ,
(pb_{m}p)a_{i_{m}} }^{( \A )}  ) ,
\]
where at the last equality sign in the sequence we used again the
Theorem 1.4 in the form presented in Remark 3.10 (and applied to the
instance ``$\{ \aunun \}$ free from $\B$''). Hence the right-hand side
of (4.17) has become:
\[
[ \mbox{coef } (1,2, \ldots ,m) ] \left(
\frac{1}{\alpha} M( \mu_{ (pb_{1}p)a_{i_{1}} , \ldots ,
(pb_{m}p)a_{i_{m}} }^{( \A )}  )  \right) ,
\]
and it is immediately verified that the left-hand side of (4.17) is exactly
the same thing. {\bf QED}

$\ $

Finally, for the proof of Application 1.10 we will need the following

$\ $

{\bf 4.10 Lemma:} Let $m$ be a positive integer, and let
$c_{1} , \ldots , c_{m} , c_{1} ' , \ldots , c_{m} '$ be random variables
in some non-commutative probability space $\ncps$. Then
\begin{equation}
[ \mbox{coef } (1, \ldots ,m) ]
( M( \mu_{ c_{1} , \ldots , c_{m} } ) \freestar
( M( \mu_{ c_{1} ' , \ldots , c_{m} ' } ))  \ =
\end{equation}
\[
= \ [ \mbox{coef } (1,2, \ldots ,2m) ]
\left( M( \mu_{ c_{1} , c_{1} ' , \ldots c_{m} , c_{m} ' } )  \freestar
Sqsum \right) ,
\]
where $Sqsum \in \Theta_{2m}$ is the quadratic polynomial
\begin{equation}
Sqsum (z_{1} , z_{2} , \ldots , z_{2m} ) \ = \
\sum_{i,j =1}^{2m} z_{i} z_{j} .
\end{equation}

$\ $

{\bf Proof} We have, by the definition,
\[
[ \mbox{coef } (1,2, \ldots ,2m) ]
\left( M( \mu_{ c_{1} , c_{1} ' , \ldots c_{m} , c_{m} ' } )  \freestar
Sqsum \right)  \ =
\]
\begin{equation}
\sum_{\sigma \in NC(2m)}
[ \mbox{coef } (1,2, \ldots , 2m) ; \sigma ]
( M( \mu_{ c_{1} , c_{1} ', \ldots , c_{m} , c_{m} ' } ) ) \cdot
[ \mbox{coef } (1,2, \ldots , 2m) ; K( \sigma ) ] (Sqsum).
\end{equation}
A partition $\sigma \in NC(2m)$ can bring a non-zero contribution to
the sum (4.22) only if $K( \sigma )$ is a pairing (otherwise
the coefficient of $Sqsum$ vanishes), so (4.22) is in fact equal to
\begin{equation}
\sum_{ \begin{array}{c}
{\scriptstyle \sigma \in NC(2m) } \\
{\scriptstyle such \ that }  \\
{\scriptstyle K( \sigma ) \in NCP(2m)}
\end{array}  }
[ \mbox{coef } (1,2, \ldots , 2m) ; \sigma ]
( M( \mu_{ c_{1} , c_{1} ', \ldots , c_{m} , c_{m} ' } ) ) .
\end{equation}
But, in view of the Proposition in 2.6, a partition $\sigma \in NC(2m)$ has
$K( \sigma ) \in NCP(2m)$ if and only if there exists a
$\rho \in NC(m)$ (uniquely determined) such that $\sigma$ is obtained by
letting $\rho$ work on $\{ 1,3, \ldots , 2m-1 \}$ and by letting
$K( \rho )$ work on $\{ 2,4, \ldots , 2m \}$. Moreover, if $\sigma \in NC(2m)$
and $\rho \in NC(m)$ are related in this way, it is immediate that
\[
[ \mbox{coef } (1,2, \ldots , 2m) ; \sigma ]
( M( \mu_{ c_{1} , c_{1} ', \ldots , c_{m} , c_{m} ' } ) ) \ = \
\]
\[
[ \mbox{coef } (1, \ldots ,m) ; \rho ]
( M( \mu_{ c_{1} , \ldots , c_{m} } ))  \cdot
[ \mbox{coef } (1, \ldots ,m) ; K( \rho ) ]
( M( \mu_{ c_{1} ' , \ldots , c_{m} ' } )) .
\]
This observation shows that the sum (4.23) can be re-written as
\[
\sum_{ \rho \in NC(m) } [ \mbox{coef } (1, \ldots ,m) ; \rho ]
( M( \mu_{ c_{1} , \ldots , c_{m} } ))  \cdot
[ \mbox{coef } (1, \ldots ,m) ; K( \rho ) ]
( M( \mu_{ c_{1} ' , \ldots , c_{m} ' } )) ;
\]
the latter quantity is just the expansion of the left-hand side of (4.20),
so we are done. {\bf QED}

$\ $

{\bf 4.11 The proof of Application 1.10} Let $\ncps$ and $\aunun ,b \in \A$
be as in the statement of (1.10). Let ${\cal C}$ denote the unital subalgebra
of $\A$ generated by $\aunun$, and put $\X = \{ ba_{1}b, \ldots , ba_{n}b \}$.
We will use the freeness criterion 4.7, in the form mentioned in Remark 4.8,
for proving that $\X$ is free from ${\cal C}$ in $\ncps$. The sufficient
condition (4.16) of 4.8 becomes here:
\begin{equation}
\varphi( c_{1} (ba_{i_{1}}b) c_{2} (ba_{i_{2}}b) \cdots
c_{m} (ba_{i_{m}}b)) \ = \ [ \mbox{coef } (1,2, \ldots ,m) ]
\left( M( \mu_{ c_{1} , \ldots , c_{m} } ) \freestar
R( \mu_{ ba_{i_{1}}b , \ldots , ba_{i_{m}}b } )  \right) ,
\end{equation}
to be verified for every $m \geq 1, \ c_{1} , \ldots , c_{m} \in {\cal C}$
and $1 \leq i_{1} , \ldots , i_{m} \leq n.$

The left-hand side of (4.24) is equal to
$[ \mbox{coef } (1,2, \ldots , 2m) ]
( M( \mu_{ c_{1}b, a_{i_{1}}b, \ldots ,c_{m}b, a_{i_{m}}b } )).$
But since $b$ is free from ${\cal C}$ (which contains
$c_{1} , \ldots , c_{m}, a_{i_{1}}, \ldots , a_{i_{m}} ),$ we have, by
Theorem 1.4 and Remark 3.10:
\[
M( \mu_{ c_{1}b, a_{i_{1}}b, \ldots ,c_{m}b, a_{i_{m}}b } ) \ = \
M( \mu_{ c_{1} , a_{i_{1}} , \ldots ,c_{m} , a_{i_{m}}  } ) \freestar
R( \mu_{ \underbrace{b, \ldots b}_{2m} } ).
\]
As it was pointed out in Lemma 4.2,
$R( \mu_{ \underbrace{b, \ldots b}_{2m} } )$ is exactly $Sqsum \circ D_{r/2}$,
where $r$ is the radius of $b$ and $Sqsum \in \Theta_{2m}$ is the series
appearing in Eqn.(4.21) of Lemma 4.10. Hence (by also using (4.3)):
\[
M( \mu_{ c_{1}b, a_{i_{1}}b, \ldots ,c_{m}b, a_{i_{m}}b } ) \ = \
\left( M( \mu_{ c_{1} , a_{i_{1}} , \ldots ,c_{m} , a_{i_{m}}  } )
\freestar Sqsum \right)  \circ D_{r/2} ,
\]
and the left-hand side of (4.24) becomes
\[
{ \left( \frac{r}{2} \right) }^{2m} \cdot
[ \mbox{coef } (1,2, \ldots , 2m) ]
\left( M( \mu_{ c_{1} , a_{i_{1}} , \ldots ,c_{m} , a_{i_{m}}  } )
\freestar Sqsum \right)  \ =
\]
\begin{equation}
\stackrel{Lemma \ 4.10}{=} \
{ \left( \frac{r}{2} \right) }^{2m} \cdot
[ \mbox{coef } (1, \ldots , m) ]
\left( M( \mu_{ c_{1} , \ldots , c_{m} } )  \freestar
M( \mu_{ a_{i_{1}} , \ldots , a_{i_{m}} } )  \right) .
\end{equation}

On the other hand, we know from Application 1.6 that
$R( \mu_{ ba_{i_{1}}b , \ldots , ba_{i_{m}}b } )$ =
$M( \mu_{ a_{i_{1}} , \ldots , a_{i_{m}} } ) \circ D_{r^{2}/4};$
consequently, we have
\[
M( \mu_{ c_{1} , \ldots , c_{m} } ) \freestar
R( \mu_{ ba_{i_{1}}b , \ldots , ba_{i_{m}}b } )  \ = \
\left( M( \mu_{ c_{1} , \ldots , c_{m} } ) \freestar
M( \mu_{ a_{i_{1}} , \ldots , a_{i_{m}} } ) \right)  \circ D_{r^{2}/4},
\]
and this makes clear that the right-hand side of (4.24) is also equal
to (4.25). {\bf QED}

$\ $

$\ $

\end{document}